\documentclass[10pt,reqno]{amsart}

\title[The meta-abelian elliptic KZB associator]{The meta-abelian elliptic KZB associator and periods of Eisenstein series}
\author{Nils Matthes}

\usepackage{amssymb,amsmath,amscd,amsthm,mathrsfs}
\usepackage{tikz}

\usepackage{color}

\usepackage{mathtools}
\mathtoolsset{showonlyrefs}

\date{}
\address{Max-Planck-Institut f\"ur Mathematik, Vivatsgasse 7, D-53111, Bonn, Germany}
\email{nilsmath@mpim-bonn.mpg.de}

\makeindex

\newif\ifnote 
\notetrue

\allowdisplaybreaks

\RequirePackage{verbatim}

\makeatletter
\newwrite\bibinl@out
\newenvironment{bibtex}[1][\jobname]{%
 \immediate\openout\bibinl@out #1.bib
 \immediate\write\bibinl@out{\@percentchar generated from `\jobname' starting line \the
\inputlineno^^J}%
 \def\verbatim@processline{\immediate\write\bibinl@out{\the\verbatim@line}}%
 \@bsphack\let\do\@makeother\dospecials\catcode`\^^M\active\verbatim@start
}%
{\immediate\closeout\bibinl@out\@esphack}
\makeatother

\theoremstyle{definition}
\newtheorem{dfn}{Definition}[section]
\newtheorem{rmk}[dfn]{Remark}

\theoremstyle{plain}
\newtheorem{prop}[dfn]{Proposition}

\newtheorem{thm}[dfn]{Theorem}
\newtheorem{lem}[dfn]{Lemma}
\newtheorem{cor}[dfn]{Corollary}

\newtheorem*{intthm}{Theorem}

\newenvironment{prf}{\begin{proof}[{\it Proof: \nopunct}]}{\end{proof}}

\numberwithin{equation}{section}

\def\dd{\mathrm{d}}

\def\bC{\mathbb C}
\def\bZ{\mathbb Z}
\def\bQ{\mathbb Q}

\def\bL{\mathbb L}
\def\bP{\mathbb P}
\def\bR{\mathbb R}

\def\cE{\mathcal E}
\def\cG{\mathcal G}

\def\cM{\mathcal M}
\def\cU{\mathcal U}

\def\fA{\mathfrak{A}}
\def\fB{\mathfrak{B}}

\def\fH{\mathfrak{H}}

\def\fp{\mathfrak{p}}

\def\fu{\mathfrak{u}}

\DeclareMathOperator{\ad}{ad}
\DeclareMathOperator{\Aut}{Aut}
\DeclareMathOperator{\gr}{gr}
\DeclareMathOperator{\Der}{Der}
\DeclareMathOperator{\End}{End}
\DeclareMathOperator{\im}{Im}

\DeclareMathOperator{\SL}{SL}

\def\sa{\mathsf{a}}
\def\sb{\mathsf{b}}
\def\sx{\mathsf{x}}
\def\sU{\mathsf{U}}
\def\sV{\mathsf{V}}
\def\sW{\mathsf{W}}
\def\sX{\mathsf{X}}
\def\sY{\mathsf{Y}}

\begin{document}

\begin{abstract}
We compute the image of Enriquez' elliptic KZB associator in the (maximal) meta-abelian quotient of the fundamental Lie algebra of a once-punctured elliptic curve. Our main result is an explicit formula for this image in terms of Eichler integrals of Eisenstein series, and is analogous to Deligne's computation of the depth one quotient of the Drinfeld associator. We also show how to retrieve Zagier's extended period polynomials of Eisenstein series, as well as the values at zero of Beilinson--Levin's elliptic polylogarithms from the meta-abelian elliptic KZB associator.
\end{abstract}

\maketitle

\section{Introduction} \label{sec:1}

This paper deals with the computation of some of the coefficients of the elliptic KZB associator defined by Enriquez \cite{Enriquez:EllAss}. In order to put things into context, we first recall the analogous picture in genus zero, due to Deligne, Drinfeld and Ihara.

Let $\fp(U):=\bL(\sx_0,\sx_1)^{\wedge}$ be the lower central series completion of the free Lie algebra in variables $\sx_0,\sx_1$, and denote by $\exp \fp(U)$ the associated pro-unipotent algebraic group. The Drinfeld associator $\Phi(\sx_0,\sx_1)$ is an element of $\exp \fp(U)_{\bR}:=\exp (\fp(U) \widehat{\otimes} \bR)$, which is constructed from the monodromy of the universal Knizhnik--Zamolodchikov (KZ) connection on $\bP^1_{\bC} \setminus \{0,1,\infty\}$ (for this reason, $\Phi$ is sometimes called KZ-associator). First introduced in \cite{Drinfeld:Gal}, the Drinfeld associator plays a pivotal role in the context of quantum groups and Grothendieck--Teichm\"uller theory.

We are interested in arithmetic properties of $\Phi(\sx_0,\sx_1)$. The following two aspects, which are in fact closely related to each other, are of particular relevance.
\begin{enumerate}
\item[(i)]
The coefficients of $\Phi(\sx_0,\sx_1)$ are expressible as $\bQ$-linear combinations of multiple zeta values
\begin{equation} \label{eqn:mzv}
\zeta(k_1,\ldots,k_n)=\sum_{m_1>\ldots>m_n>0}\frac{1}{m_1^{k_1}\ldots m_n^{k_n}}, \quad k_1 \geq 2, \, k_2,\ldots,k_n \geq 1,
\end{equation}
which are generalizations of the special values of the Riemann zeta function at positive integers. These numbers have (at least conjecturally) a rich algebraic structure \cite{Goncharov:MTM,IKZ}.
\item[(ii)]
The Lie algebra $\fp(U)$ is the de Rham realization of an element of the category $\mathsf{MTM}$ of mixed Tate motives over $\bZ$ (\cite{DG}, \S 5). As a consequence, the unipotent fundamental group $\cU_{\mathsf{MTM}}$ of $\mathsf{MTM}$ acts on $\exp\fp(U)$ (Ihara action), and in particular on $\Phi(\sx_0,\sx_1)$.\footnote{In this context, $\Phi(\sx_0,\sx_1)$ is usually denoted $dch$ (for `droit chemin').} The Deligne--Ihara conjecture (proved by Brown in \cite{Brown:MTM}) states that this action is faithful, thus elements of $\cU_{\mathsf{MTM}}$ are completely determined by their action on $\Phi(\sx_0,\sx_1)$, which can be computed very explicitly \cite{Brown:Decomposition}.
\end{enumerate}
For both (i) and (ii), the archetypal result is due to Deligne (\cite{Deligne:P1}, \S 19), who inspired by unpublished work of Wojtkowiak essentially showed that
\begin{equation} \label{eqn:depthone}
\log(\Phi(\sx_0,\sx_1)) \equiv -\sum_{k=2}^{\infty}\zeta(k)\ad^{k-1}(\sx_0)(\sx_1) \mod [D^1\fp(U),D^1\fp(U)],
\end{equation}
where $D^1\fp(U) \subset \fp(U)$ denotes the ideal generated by $\sx_1$. On the one hand, this exhibits the Riemann zeta values $\zeta(k)$ as coefficients of $\log(\Phi(\sx_0,\sx_1))$. On the other hand, since $\zeta(k) \neq 0$, one deduces from \eqref{eqn:depthone} that the generators $\exp(\sigma_{2n+1})$ of $\cU_{\mathsf{MTM}}$ act non-trivially on $\exp \fp(U)$ (\cite{DG}, \S 6.8), which was a first step towards establishing the Deligne--Ihara conjecture.

In this paper, we consider an elliptic analog of the above situation. Let $\fH$ be the Poincar\'e upper half-plane, and consider for $\tau \in \fH$ the once-punctured, complex elliptic curve $E_{\tau}^{\times}:=\bC/(\bZ+\bZ\tau) \setminus \{0\}$. Following Hain--Matsumoto \cite{HM}, we denote its de Rham fundamental group by $\fp(E_{\tau}^{\times}) \cong \bL(\sa,\sb)^{\wedge}$. In \cite{Enriquez:EllAss}, Enriquez constructs the elliptic KZB associator $(A(\tau),B(\tau)) \in \exp \fp(E_{\tau}^{\times})_{\bC} \times \exp \fp(E_{\tau}^{\times})_{\bC}$ from the monodromy of the universal elliptic Knizhnik--Zamolodchikov--Bernard (KZB) connection \cite{CEE:KZB,LR}. It is an elliptic version of the Drinfeld associator, and the analogs of (i) and (ii) above are the following.
\begin{enumerate}
\item[(i)]
The coefficients of the elliptic KZB associator are the elliptic multiple zeta values, first introduced in \cite{Enriquez:Emzv} and studied in more detail in \cite{BMS,LMS,Matthes:Thesis,Matthes:Edzv}.
They are closely related to both multiple zeta values and to iterated integrals of Eisenstein series \cite{Brown:MMV,Manin:Iterated}.
\item[(ii)] The Lie algebra $\fp(E_{\tau}^{\times})$, viewed as a local system over the moduli space $\cM_{1,\overrightarrow{1}}$ of elliptic curves with a non-zero tangent vector at the origin, is the de Rham realization of an element of the category $\mathsf{MEM}_{\overrightarrow{1}}$ of universal mixed elliptic motives (over $\cM_{1,\overrightarrow{1}}$). This category can be seen as an elliptic enhancement of the category of mixed Tate motives over $\bZ$. The corresponding Galois group $\cG_{\mathsf{MEM}_{\overrightarrow{1}}}$ acts on $\fp(E_{\tau}^{\times})$ \cite{HM}, and therefore also on the elliptic KZB associator. In analogy to the Deligne--Ihara conjecture, it is asked in \cite{HM}, \S 24.2 whether the action of $\cG_{\mathsf{MEM}_{\overrightarrow{1}}}$ on $\fp(E_{\tau}^{\times})$ is faithful.
\end{enumerate}
The main goal of this article is to establish an analog of \eqref{eqn:depthone} for the elliptic KZB associator, i.e. the explicit computation of the images of the formal logarithms $\fA(\tau):=\log(A(\tau))$ and $\fB(\tau):=\log(B(\tau))$ in a certain quotient of $\fp(E_{\tau}^{\times})_{\bC}$. More precisely, let $D^1\fp(E_{\tau}^{\times}) \subset \fp(E_{\tau}^{\times}) \cong \bL(\sa,\sb)^{\wedge}$ be the commutator. Taking its lower central series defines a filtration $D^{\bullet}\fp(E_{\tau}^{\times})$, the elliptic depth filtration (\cite{HM}, \S 27). In particular, $D^2\fp(E_{\tau}^{\times})$ is the double commutator, and our goal is to compute the images $\fA(\tau)^{\rm met-ab}$ and $\fB(\tau)^{\rm met-ab}$ of the elliptic KZB associator in the meta-abelian quotient 
\begin{equation}
\fp(E_{\tau}^{\times})^{\rm met-ab}_{\bC}:=\fp(E_{\tau}^{\times})_{\bC}/D^2\fp(E_{\tau}^{\times})_{\bC} \cong (\bC\cdot \sa\oplus \bC\cdot \sb) \oplus \bC[\![U,V]\!],
\end{equation}
where $U^kV^l:=\ad^k(\sa)\ad^l(\sb)([\sa,\sb])$. Our main result can then be stated as follows.
\begin{intthm}[Theorem \ref{thm:arithgeo} below]
Let $\overline{\sU}:=\frac{\sU}{2\pi i}$ and $\sW:=\overline{\sU}+\tau \sV$. We have
\begin{align}
\fA(\tau)^{\rm met-ab}&=2\pi i\sb+\exp\left(\tau\frac{\partial}{\partial \overline{\sU}}\sV \right)\fA^{(1)}_{\infty}-2\pi i\sV\sum_{k=1}^{\infty}\frac{2}{(2k-2)!}\int_{\tau}^{\overrightarrow{1}_{\infty}}\underline{G}_{2k},
\end{align}
and
\begin{align}
\fB(\tau)^{\rm met-ab}&=\sa+2\pi i\tau \sb+\exp\left(\tau\frac{\partial}{\partial \overline{\sU}}\sV \right)\fB^{(1)}_{\infty}-2\pi i\sW\sum_{k=1}^{\infty}\frac{2}{(2k-2)!}\int_{\tau}^{\overrightarrow{1}_{\infty}}\underline{G}_{2k}.
\end{align}
Here, $\int_{\tau}^{\overrightarrow{1}_{\infty}}\underline{G}_{2k}:=(2\pi i)^{2k-1}\int_{\tau}^{\overrightarrow{1}_{\infty}}G_{2k}(z)(\sW-z\sV)^{2k-2}\dd z$ is the regularized Eichler integral of $G_{2k}$ (\cite{Brown:MMV}, \S 4), and the series $\fA_{\infty}^{(1)}$, $\fB_{\infty}^{(1)}$ are given by
\begin{align}
\fA_{\infty}^{(1)}&=2\pi i\left(c(\sU)-\frac{(2\pi i)}{4} \sV+\sum_{n \geq 3, {\rm odd}}\zeta(n)\sV^n\right),\\
\fB_{\infty}^{(1)}&=-2\pi i\left(c(2\pi i\sV)-\sU c(\sU)c(2\pi i\sV)\right)+\sum_{n \geq 3, \,{\rm odd}}\zeta(n)\sU\sV^{n-1},
\end{align}
where $c(x)=\frac{1}{e^x-1}+\frac 12-\frac 1x=\sum_{k=2}^{\infty}\frac{B_k}{k!}x^{k-1}$.
\end{intthm}
Similar considerations have been made by Hain to prove that the generators $\exp(\mathbf{e}_{2k})$ of the geometric fundamental group $\cG^{\rm geom}_{\mathsf{MEM}_{\overrightarrow{1}}}$ act non-trivially on $\fp(E_{\tau}^{\times})$ (\cite{Hain:HodgeDeRham}, Theorem 15.7).
Moreover, our theorem gives a closed expression of elliptic multiple zeta values of depth one explicitly in terms of Riemann zeta values and Eichler integrals of Eisenstein series.


The proof of Theorem \ref{thm:arithgeo} uses a result of Enriquez \cite{Enriquez:EllAss} to the effect that 
\begin{equation} \label{eqn:crucial}
\fA(\tau)=g(\tau)(\fA_{\infty}), \quad \fB(\tau)=g(\tau)(\fB_{\infty}),
\end{equation}
for certain explicit elements $\fA_{\infty},\fB_{\infty} \in \fp(E_{\tau}^{\times})_{\bC}$ and an automorphism $g(\tau) \in \Aut(\exp(\fp(E_{\tau}^{\times})_{\bC}))$. Then, we separately compute the images of $\fA_{\infty}$ and $\fB_{\infty}$ in $\fp(E_{\tau}^{\times})_{\bC}^{\rm met-ab}$ and of $g(\tau)$ in $\Aut(\exp(\fp(E_{\tau}^{\times})_{\bC}^{\rm met-ab}))$, and from this, we are able to deduce Theorem \ref{thm:arithgeo}.

The series $\fA_{\infty}$ and $\fB_{\infty}$ are arithmetic: they can be expressed in terms of the Drinfeld associator and therefore come from genus zero. On the other hand, the automorphism $g(\tau)$ is geometric: it describes the action of $\cG^{\rm geom}_{\mathsf{MEM}_{\overrightarrow{1}}}$ on $\exp \fp(E_{\tau}^{\times})$. As a byproduct of our proof, we see that already their images in the meta-abelian quotient are interesting objects in their own right. Namely, the automorphism $g(\tau)^{\rm met-ab}$ is essentially the generating series of the special values of elliptic polylogarithms at the zero section of the elliptic curve \cite{BeiLev,Levin:Compositio} (cf. Theorem \ref{thm:geometric} and Corollary \ref{cor:geometric}), while $\fA_{\infty}^{\rm met-ab}$, $\fB_{\infty}^{\rm met-ab}$ turn out to be generating series of the extended period polynomials of Eisenstein series \cite{Zagier:Periods} (cf. Theorem \ref{thm:arithmetic} and Corollary \ref{cor:arithmetic}). 

Finally, we note that Nakamura \cite{Nakamura:Galoisrep,Nakamura:Eisenrevisited} has studied an $\ell$-adic analog of the meta-abelian image of the elliptic KZB associator (called ``universal power series for Dedekind sums''), which is a genus one analog of Ihara's universal power series for Jacobi sums \cite{Ihara:Annals}. It would be very interesting to compare his results to ours.

The plan of the paper is as follows. In Sections \ref{sec:2} and \ref{sec:3}, we collect some background in order to make the paper self-contained. Then, in Section \ref{sec:4}, we recall the definition of the elliptic KZB associator \cite{Enriquez:EllAss}, but from the point of view of the mixed Hodge structure on the unipotent fundamental group of $E_{\tau}^{\times}$ \cite{BL:MEP}. Finally, in Section \ref{sec:5}, the main results of this paper are proved.

{\bf Acknowledgments:} Very many thanks to B. Enriquez and H. Nakamura for very inspiring discussions at the conference ``GRT, MZVs and associators'' in Les Diablerets in 2015, which formed the starting point of this project. Thanks are also due to A. Alekseev for the invitation to that conference. Also, many thanks to F. Brown, B. Enriquez, H. Furusho and F. Zerbini for helpful comments on an earlier version of this paper. This paper was written while the author was a Ph.D. student at Universit\"at Hamburg under the supervision of U. K\"uhn.

\section{Preliminaries} \label{sec:2}

\subsection{Notation and conventions}

We start by introducing some general notation, to be used throughout the text.

We denote by $\fH:=\{ z \in \bC \, \vert \, \im(z)>0 \}$ the upper half-plane, with canonical coordinate $\tau$. For $\tau \in \fH$, we let $E_{\tau}^{\times}:=\bC/(\bZ+\bZ\tau) \setminus \{0\}$ be the associated once-punctured complex elliptic curve.

For any finite set $\{\sx_1,\ldots,\sx_n\}$ and a field $K$, we denote by $\bL(\sx_1,\ldots,\sx_n)_K$ the free Lie algebra on $X$ over $K$ (we omit $K$ if $K=\bQ$), and by $\bL(\sx_1,\ldots,\sx_n)^{\wedge}_K$ the completion for its lower central series. It is a topological Lie algebra over $K$, whose topology is induced from the lower central series. Its topological universal enveloping algebra is given by $K\langle\!\langle \sx_1,\ldots,\sx_n\rangle\!\rangle$, the $K$-algebra of formal power series in the non-commuting variables $\sx_1,\ldots,\sx_n$, and the exponential map $\exp: \bL(\sx_1,\ldots,\sx_n)^{\wedge}_K \rightarrow K\langle\!\langle \sx_1,\ldots,\sx_n\rangle\!\rangle$ defines an isomorphism onto the subspace of $K\langle\!\langle \sx_1,\ldots,\sx_n\rangle\!\rangle$ of group-like elements, denoted by $\exp \bL(\sx_1,\ldots,\sx_n)^{\wedge}_K$. For more background, we refer to \cite{Reutenauer,Serre:Lie}.

\subsection{Derivations on the fundamental Lie algebra of a once-punctured elliptic curve}

Following \cite{HM}, we will denote by $\fp(E^{\times}_{\tau})$ the (de Rham) fundamental Lie algebra of the once-punctured elliptic curve $E_{\tau}^{\times}$. With notation as above, one has
\begin{equation}
\fp(E^{\times}_{\tau})\cong\bL(\sa,\sb)^{\wedge}
\end{equation}
where the generators $\sa,\sb$ correspond to the natural homology cycles on $E_{\tau}^{\times}$.

We will need to consider a special family of derivations on $\fp(E_{\tau}^{\times})$. Denote by $\Der^0(\fp(E_{\tau}^{\times}))$ the Lie algebra of continuous derivations $D$, which satisfy $D([\sa,\sb])=0$ and such that $D(\sb)$ has no linear term in $\sa$. From these two conditions, it follows easily that every $D \in \Der^0(\fp(E_{\tau}^{\times}))$ is uniquely determined by its value on $\sa$.
\begin{dfn}[Tsunogai] \label{dfn:derivations}
For every $k \geq 0$, define $\varepsilon_{2k} \in \Der^0(\fp(E_{\tau}^{\times}))$ by its value on $\sa$:
\begin{equation} \label{eqn:valueonx}
\varepsilon_{2k}(\sa)=\begin{cases}-\sb & k=0\\ \frac{2}{(2k-2)!}\ad^{2k}(\sa)(\sb) & k>0.\end{cases}
\end{equation}
We also let $\fu \subset \Der^0(\fp(E_{\tau}^{\times}))$ be the Lie subalgebra generated by the $\varepsilon_{2k}$.
\end{dfn}
The derivations $\varepsilon_{2k}$ have first been introduced by Tsunogai (\cite{Tsunogai:Derivations}, \S 3) in the context of Galois actions on fundamental groups of punctured elliptic curves. They also play an important role in the theory of universal mixed elliptic motives, as the relative unipotent completion of $\SL_2(\bZ)$ acts on $\fp(E_{\tau}^{\times})$ through them (\cite{HM}, \S 20).
\begin{rmk}
The value of $\varepsilon_{2k}$ on $\sb$ is given by
\begin{equation}
\varepsilon_{2k}(\sb)=\frac{2}{(2k-2)!}\sum_{0 \leq j <k}(-1)^j[\ad^j(\sa)(\sb),\ad^{2k-1-j}(\sa)(\sb)].
\end{equation}
In particular, $\varepsilon_0(\sb)=0$.
\end{rmk}

\subsection{Eichler integrals of Eisenstein series} \label{ssec:2.2}

Consider the Hecke-normalized Eisenstein series for $\SL_2(\bZ)$ of weight $2k$:
\begin{equation} \label{eqn:Eis}
G_{2k}(q):=\begin{cases}-\frac{B_{2k}}{4k}+\sum_{n=1}^{\infty}\left( \sum_{d\vert n}d^{2k-1} \right)q^n & k \geq 1 \\ -1 & k=0,\end{cases}
\end{equation}
where $B_{2k}$ denotes the $2k$-th Bernoulli number and $q=e^{2\pi i\tau}$. Extending earlier work of Manin \cite{Manin:Iterated}, Brown \cite{Brown:MMV} introduced (regularized) iterated integrals of \eqref{eqn:Eis} (or \textit{iterated Eisenstein integrals} for short) 
\begin{equation} \label{eqn:IterEis}
\cG(2k_1,\ldots,2k_n;\tau):=\int_{\tau}^{\overrightarrow{1}_{\infty}}G_{2k_1}(\tau_1)\ldots  G_{2k_n}(\tau_n) \dd\tau_1\ldots\dd\tau_n,
\end{equation}
where $\overrightarrow{1}_{\infty}$ denotes the tangential base point $1$ at $i\infty$. We refer to \cite{Brown:MMV}, \S 4, for the general definition, and only note the special case
\begin{align}
\cG(\{0\}_n,2k;\tau)&=(-1)^n\idotsint\limits_{\tau \leq \tau_1\leq \ldots \tau_{n+1} \leq i\infty}G_{2k}(\tau_{n+1})-a_0(G_{2k})\dd \tau_1\ldots\dd\tau_{n+1}\notag\\
&-a_0(G_{2k})\frac{\tau^{n+1}}{(n+1)!},\label{eqn:IEIspecial1}
\end{align}
where $\{0\}_n$ denotes an $n$-tuple of zeros, and $a_0(G_{2k})=-\frac{B_{2k}}{4k}$ is the constant term in the Fourier expansion \eqref{eqn:Eis} of $G_{2k}$. From the shuffle product formula for (regularized) iterated integrals (\cite{Brown:MMV}, Proposition 4.7), we further deduce
\begin{align} 
\cG(\{0\}_{n-1},2k,0;\tau)&=\cG(0;\tau)\cG(\{0\}_{n-1},2k;\tau)-n\cG(\{0\}_n,2k;\tau) \label{eqn:IEIspecial2}.
\end{align}
Both $\cG(\{0\}_n,2k;\tau)$ and $\cG(\{0\}_{n-1},2k,0;\tau)$ can be expressed in terms of generalized Eichler integrals
\begin{equation} \label{eqn:EichlerIntegral}
I_n(G_{2k};\tau):=\int_{\tau}^{i\infty}\Big[G_{2k}(z)-a_0(G_{2k})\Big](\tau-z)^n\dd z-\int_0^{\tau}a_0(G_{2k})(\tau-z)^n\dd z,
\end{equation}
with the classical Eichler integral of $G_{2k}$ being the special case $n=2k-2$ and $k \geq 2$ (cf. e.g. \cite{Zagier:Traces}, \S 1).
\begin{prop} \label{prop:IterEisEichler}
We have
\begin{align}
\cG(\{0\}_n;\tau)&=\frac{\tau^n}{n!} \label{eqn:1}\\
\cG(\{0\}_n,2k;\tau)&=\frac{1}{n!}I_n(G_{2k};\tau) \label{eqn:2},
\end{align}
and for $k,n\geq 1$:
\begin{equation} \label{eqn:3}	\cG(\{0\}_{n-1},2k,0;\tau)=\frac{1}{(n-1)!}\left(\tau I_{n-1}(G_{2k};\tau)-I_n(G_{2k};\tau)\right).
\end{equation}
\end{prop}
\begin{prf}
The first equality is immediate from the definition \eqref{eqn:IEIspecial1}. The second equality \eqref{eqn:2} is trivial for $n=0$, and the general case is easy to prove from \eqref{eqn:1} by induction on $n$. Finally, \eqref{eqn:3} follows directly from \eqref{eqn:1}, \eqref{eqn:2} and the definition \eqref{eqn:IEIspecial2}.
\end{prf}

\subsection{The elliptic KZB connection and the associated transport map} \label{ssec:2.3}


We recall the definition of the elliptic KZB (Knizhnik--Zamolodchikov--Bernard) connection $\nabla_{\rm KZB}$ on $E_{\tau}^{\times}$, whose monodromy will give rise to the elliptic KZB associator. Originally, $\nabla_{\rm KZB}$ was defined as a meromorphic connection on $\bC$ (cf. \cite{CEE:KZB,Hain:KZB,LR}). Here, we will instead follow \cite{BL:MEP}, which consider a certain $C^{\infty}$-trivialization of $\nabla_{\rm KZB}$, which is defined on the quotient $\bC/(\bZ+\bZ\tau) \setminus \{0\}$.

Let $\xi=r\tau+s$ be the canonical coordinate on $E_{\tau}^{\times}$, with $(r,s) \in \bR^2 \setminus \bZ^2$. Also, let
\begin{equation}
\theta_{\tau}(\xi)=\sum_{n \in \bZ}(-1)^nq^{\frac{1}{2}(n+\frac 12)^2}e^{(n+\frac 12)\xi}, \quad q=e^{2\pi i\tau},
\end{equation}
be the classical Jacobi theta function.
\begin{dfn}[Brown--Levin,Calaque--Enriquez--Etingof,Levin--Racinet]
Define a connection $\nabla_{\rm KZB}$ on the trivial bundle\footnote{Note that the normalization of the variables $\sa,\sb$ differs from \cite{BL:MEP}, Example 5.3.1, by $\sa=-2\pi i\mathsf{x}_0$ and $\sb=-(2\pi i)^{-1}\mathsf{x}_1$. Our conventions are compatible with \cite{Hain:KZB}, \S 11.1.}
\begin{equation}
E_{\tau}^{\times} \times \bC\langle\!\langle \sa,\sb\rangle\!\rangle \rightarrow E_{\tau}^{\times}
\end{equation}
by setting $\nabla_{\rm KZB}(f):=\dd f-\omega_{\rm KZB}\cdot f$ for a local section $f$, where
\begin{equation} \label{eqn:KZBform}
\omega_{\rm KZB}=\dd r\cdot \sa+2\pi i\ad(\sa)e^{r\ad(\sa)}F_{\tau}(2\pi i\xi,\ad(\sa))(\sb)\dd \xi,
\end{equation}
where
\begin{equation}
F_{\tau}(\xi,\eta):=\frac{\theta'_{\tau}(0)\theta_{\tau}(\xi+\eta)}{\theta_{\tau}(\xi)\theta_{\tau}(\eta)}.
\end{equation}
\end{dfn}

\begin{prop} \label{prop:connprop}
The connection $\nabla_{\rm KZB}$ satisfies the following properties.
\begin{itemize}
\item[\rm (i)]
We have $\nabla_{\rm KZB}^2=0$; in other words, $\nabla_{\rm KZB}$ is integrable.
\item[\rm (ii)]
The connection $\nabla_{\rm KZB}$ has a simple pole at $\xi=0$ with residue
\begin{equation}
\operatorname{Res}_0(\nabla_{\rm KZB})=[\sa,\sb].
\end{equation}
\end{itemize}
\end{prop}
\begin{prf}
\begin{itemize}
\item[\rm (i)]
The condition $\nabla_{\rm KZB}^2=0$ is equivalent to
\begin{equation}
\dd\omega_{\rm KZB}-\omega_{\rm KZB}\wedge \omega_{\rm KZB}=0,
\end{equation}
which in turn follows from a direct computation:
\begin{align}
\dd\omega_{\rm KZB}&=2\pi i\dd r \cdot \ad(\sa) \wedge \ad(\sa)e^{r\ad(\sa)}F_{\tau}(2\pi i\xi,\ad(\sa))(\sb)\dd \xi\\
&=\omega_{\rm KZB}\wedge \omega_{\rm KZB}.
\end{align}
\item[\rm (ii)]
The residue of the connection $\nabla_{\rm KZB}$ is just the residue of the one-form $\omega_{\rm KZB}$. But the computation of the latter is easy from the definition, using the fact that the residue of $2\pi iF_{\tau}(2\pi i\xi,\eta)$ at $\xi=0$ is equal to one (cf. \cite{Hain:KZB}, eqn.(8)).
\end{itemize}
\end{prf}
Now for any two base points $\rho_1,\rho_2$, let $\pi_1(E_{\tau}^{\times};\rho_2,\rho_1)$ be the fundamental torsor of paths from $\rho_1$ to $\rho_2$. The integrability of $\nabla_{\rm KZB}$ implies that the transport function
\begin{align}
T^{\rm KZB}_{\rho_2,\rho_1}: \pi_1(E_{\tau}^{\times};\rho_2,\rho_1) &\rightarrow \bC\langle\!\langle \sa,\sb\rangle\!\rangle\\
\gamma &\mapsto \sum_{k=0}^{\infty}\int_{\gamma}\omega^k_{\rm KZB},
\end{align}
is well-defined, where $\int_{\gamma}\omega_{\rm KZB}^k$ denotes the iterated integral in the sense of Chen \cite{Chen:PathIntegrals}
\begin{equation}
\int_{\gamma}\omega_{\rm KZB}^k:=\int\limits_{1\geq t_1\geq \ldots \geq t_k \geq 1}\gamma^*(\omega_{\rm KZB})(t_1)\ldots \gamma^*(\omega_{\rm KZB})(t_k).
\end{equation}
In other words, $\int_{\gamma}\omega_{\rm KZB}^k$ depends only on the homotopy class of $\gamma$.

Rather than choosing points $\rho_1,\rho_2 \in E_{\tau}^{\times}$, which is not canonical, we work with tangential base points, in the sense of \cite{Deligne:P1}, \S 15, at the puncture $0$. Since $\nabla_{\rm KZB}$ has only a simple pole at $\xi=0$, one can extend the definition of the transport function to the case of tangential base points as in \cite{Deligne:P1}, Proposition 15.45. More precisely, for any two non-zero tangent vectors $\overrightarrow{v}_0=\lambda \frac{\partial}{\partial \xi}$ and $\overrightarrow{w}_0=\mu\frac{\partial}{\partial \xi}$ at $0$, there is a well-defined function
\begin{align}
T^{\rm KZB}_{\overrightarrow{w}_0,\overrightarrow{v}_0}: \pi_1(E_{\tau}^{\times};\overrightarrow{w}_0,\overrightarrow{v}_0) &\rightarrow \bC\langle\!\langle \sa,\sb\rangle\!\rangle,
\end{align}
given by
\begin{align}
T^{\rm KZB}_{\overrightarrow{w}_0,\overrightarrow{v}_0}(\gamma)=\lim_{t\to 0}e^{\log(\mu^{-1}t)\operatorname{Res}_0(\nabla_{\rm KZB})}\Bigg[\sum_{k=0}^{\infty}\int_{\gamma_t^{1-t}}\omega^k_{\rm KZB}\Bigg] e^{-\log(\lambda^{-1}t)\operatorname{Res}_0(\nabla_{\rm KZB})},
\end{align}
where $\operatorname{Res}_0(\nabla_{\rm KZB})=[\sa,\sb]$ is the residue of the connection at $\xi=0$ (cf. Proposition \ref{prop:connprop}.(i)), $\gamma_t^{1-t}$ denotes the restriction of $\gamma$ to the interval $[t,1-t]$ (for $0<t<\frac 12$) and the branches of the logarithms are determined by the path $\gamma$. For arithmetic applications, it will be important that the tangent vectors are integral on the Tate curve $\bC^{\times}/q^{\bZ}$ and moreover non-zero modulo every prime number $p$, which fixes them uniquely (up to a sign): $\overrightarrow{v}_0=\pm \frac{\partial}{\partial z}=\pm (2\pi i)^{-1}\frac{\partial}{\partial \xi}$, where $z=e^{2\pi i\xi}$.

\section{The elliptic depth filtration} \label{sec:3}

We recall the definition of the elliptic depth filtration on the fundamental Lie algebra of $E_{\tau}^{\times}$ (cf. \cite{HM}, \S 27). This filtration is the elliptic analog of the depth filtration on the fundamental Lie algebra of $\bP^1 \setminus \{0,1,\infty\}$ (\cite{DG}, \S 6 or \cite{Brown:Depth}, \S 4).
\subsection{The elliptic depth filtration}
Consider the canonical embedding
\begin{equation} \label{eqn:embedding}
E_{\tau}^{\times} \hookrightarrow E_{\tau}
\end{equation}
of the once-punctured elliptic curve $E_{\tau}^{\times}$ into the (complete) elliptic curve $E_{\tau}$. On fundamental Lie algebras, it induces the abelianization map
\begin{equation}
\pi: \fp(E_{\tau}^{\times}) \rightarrow \fp(E_{\tau}^{\times})^{\rm ab} \cong \fp(E_{\tau}).
\end{equation}
\begin{dfn}[Hain--Matsumoto]
The \textit{elliptic depth filtration} $D^{\bullet}\fp(E_{\tau}^{\times})$ is the descending filtration on $\fp(E_{\tau}^{\times})$, defined by
\begin{equation}
D^n\fp(E_{\tau}^{\times})=\begin{cases}\fp(E_{\tau}^{\times}) & n=0\\ \ker(\pi) & n=1 \\ [D^1\fp(E_{\tau}^{\times}),D^{n-1}\fp(E_{\tau}^{\times})] & n \geq 2 \end{cases}.
\end{equation}
Also, let $\gr^{\bullet}_D\fp(E_{\tau}^{\times})$ be the associated graded Lie algebra.
\end{dfn}
It is clear from the definition that the elliptic depth filtration is the lower central series on the commutator of $\fp(E_{\tau}^{\times})$. Therefore, the quotient Lie algebra
\begin{equation}
\fp(E_{\tau}^{\times})^{\rm met-ab}:=\fp(E_{\tau}^{\times})/D^2\fp(E_{\tau}^{\times})
\end{equation}
is the \textit{(maximal) meta-abelian quotient} of $\fp(E_{\tau}^{\times})$.

The following proposition is well-known. 
\begin{prop} \label{prop:metablie}
We have isomorphisms of (abelian) Lie algebras
\begin{equation} \label{eqn:isom1}
\gr^0_D\fp(E_{\tau}^{\times}) \cong \bQ \sa\oplus \bQ \sb
\end{equation}
and
\begin{align} 
\gr^1_D\fp(E_{\tau}^{\times}) &\stackrel{\cong}\longrightarrow \bQ[\![\sU,\sV]\!] \notag\\
\ad^k(\sa)\ad^l(\sb)([\sa,\sb]) &\mapsto \sU^k\sV^l. \label{eqn:isom2}
\end{align}
Moreover, 
\begin{equation}
\fp(E_{\tau}^{\times})^{\rm met-ab} \cong \gr^0_D\fp(E_{\tau}^{\times}) \ltimes \gr^1_D\fp(E_{\tau}^{\times})
\end{equation}
as Lie algebras, where $\bQ \sa\oplus \bQ \sb$ acts on $\gr^1_D\fp(E_{\tau}^{\times})\cong \bQ[\![\sU,\sV]\!]$ by the adjoint action.
\end{prop}
\begin{prf}
The first isomorphism is clear, since the right hand side of \eqref{eqn:isom1} is just the abelianization of $\fp(E_{\tau}^{\times})$. It follows from the Jacobi identity that every element of $\gr^1_D\fp(E_{\tau}^{\times})$ is a series in the elements $\ad^k(\sa)\ad^l(\sb)([\sa,\sb])$, and then the isomorphism \eqref{eqn:isom2} is a consequence of the universal property of free Lie algebras. Finally, the last statement of the proposition follows from the fact that the adjoint action splits the short exact sequence of Lie algebras
\begin{equation}
0 \longrightarrow \gr^1_D\fp(E_{\tau}^{\times}) \longrightarrow \fp(E_{\tau}^{\times})/D^2\fp(E_{\tau}^{\times}) \longrightarrow \gr^0_D\fp(E_{\tau}^{\times}) \longrightarrow 0.
\end{equation}
\end{prf}
\begin{rmk} \label{rmk:elldepth}
The relation between the elliptic depth filtration and the depth filtration on the fundamental Lie algebra of $\bP^1 \setminus \{0,1,\infty\}$ can be explained as follows. First, recall (cf. \cite{DG}, \S 5) that the (de Rham) fundamental Lie algebra $\fp(U)$ of $U:=\bP^1 \setminus \{0,1,\infty\}$ is isomorphic to $\bL(\sx_0,\sx_1)^{\wedge}$. The depth filtration $D^n\fp(U)$ on $\fp(U)$ is then the lower central series on the kernel of the natural map between fundamental Lie algebras
\begin{align}
\fp(U) &\rightarrow \bL(\sx_0)^{\wedge} \cong \bQ \sx_0\\
\sx_i &\mapsto \delta_{i,0}\sx_0,
\end{align}
which is induced from the embedding $\bP^1 \setminus \{0,1,\infty\} \hookrightarrow \bP^1 \setminus \{0,\infty\}$ (cf. \cite{Brown:Depth,DG}). Interpreting $\bP^1 \setminus \{0,1,\infty\}$ as the fiber over $q=0$ of the universal once-punctured Tate curve $(\bC^{\times}/q^{\bZ}) \setminus \{1\}$, one obtains a morphism of Lie algebras \cite{Brown:Depth3,Enriquez:EllAss,Hain:KZB}
\begin{align}
\iota:\fp(U) &\rightarrow \fp(E_{\tau}^{\times}) \label{eqn:Hainmorphism}\\
\sx_0 &\mapsto \frac{\ad(\sa)}{e^{\ad(\sa)}-1}(\sb)=\sum_{k=0}^{\infty}\frac{B_k}{k!}\ad^k(\sa)(\sb)\\
\sx_1 &\mapsto [\sa,\sb],
\end{align}
which clearly respects the depth filtrations on both sides, i.e.
\begin{equation}
\iota(D^n\fp(U))=\iota(\fp(U))\cap D^n\fp(E_{\tau}^{\times}), \quad \mbox{for all $n \geq 0$}.
\end{equation}
For more details, see \cite{HM}, \S 27.
\end{rmk}

\subsection{Action of special derivations in depths zero and one} \label{ssec:3.2}

We now compute the action of the derivations $\varepsilon_{2k}$ on the meta-abelian quotient $\fp(E_{\tau}^{\times})^{\rm met-ab}$.
\begin{prop} \label{prop:metabaction}
\begin{enumerate}
\item[\rm (i)] The derivation $\varepsilon_0$ acts on $\gr^0_D\fp(E_{\tau}^{\times}) \cong \bQ \sa\oplus \bQ \sb$ as the linear map $\left(\begin{smallmatrix}0&-1\\0&0\end{smallmatrix}\right)$, and on $\gr^1_D\fp(E_{\tau}^{\times}) \cong \bQ[\![\sU,\sV]\!]$ as the derivation $
-\sV\frac{\partial }{\partial \sU}.
$
\item[\rm (ii)]
The derivations $\varepsilon_{2k}$, for $k>0$, act trivially on $\gr^i_D\fp(E_{\tau}^{\times})$, for every $i\geq 0$.
\item[\rm (iii)]
Let $\underline{2k}=(2k_1,\ldots,2k_n)$ be a multi-index, where $k_i \geq 0$. Then $\varepsilon_{\underline{2k}}=\varepsilon_{2k_1} \circ\ldots\circ \varepsilon_{2k_n}$ acts non-trivially on $\fp(E_{\tau}^{\times})^{\rm met-ab}\cong\gr^0_D\fp(E_{\tau}^{\times}) \ltimes \gr^1_D\fp(E_{\tau}^{\times})$, only if either $\underline{2k}=(0,\ldots,0,2k_n)$ or $\underline{2k}=(0,\ldots,0,2k_{n-1},0)$.
\end{enumerate}
\end{prop}
\begin{prf}
The action of $\varepsilon_0$ on $\gr^0_D\fp(E_{\tau}^{\times})$ is clear from the definition (cf. Definition \ref{dfn:derivations}). For the action on $\gr^1_D\fp(E_{\tau}^{\times})$, by the Jacobi identity, the linear operators $\ad(\sa),\ad(\sb) \in \End(\gr^1_D\fp(E_{\tau}^{\times}))$ commute with each other. Consequently, we have
\begin{align}
\varepsilon_0(\ad^k(\sa)\ad^l(\sb)([\sa,\sb])) &\equiv \sum_{i=0}^{k-1} -\ad^i(\sa)\ad(\sb)\ad^{k-1-i}(\sa)\ad^l(\sb)([\sa,\sb]) \notag\\
&\equiv -k\ad^{k-1}(\sa)\ad^{l+1}(\sb)([\sa,\sb]) \mod D^2\fp(E_{\tau}^{\times}).
\end{align}
Therefore, under the isomorphism $\gr^1_D\fp(E_{\tau}^{\times}) \cong \bQ[\![\sU,\sV]\!]$ of Proposition \ref{prop:metablie}, the derivation $\varepsilon_0$ corresponds to $
-\sV\frac{\partial }{\partial \sU}$. As for (ii), the triviality of $\varepsilon_{2k}$, for $k>0$, on $\gr^0_D\fp(E_{\tau}^{\times})$ is clear from Definition \ref{dfn:derivations}, and triviality on $\gr^i_D\fp(E_{\tau}^{\times})$ follows by induction on $i$. Finally, (iii) follows easily from (i) and (ii).
\end{prf}

\section{The elliptic KZB associator} \label{sec:4}

In this section, we define Enriquez's elliptic KZB associator \cite{Enriquez:EllAss}, which is an elliptic analogue of the Drinfeld associator \cite{Drinfeld:Gal}. Our approach differs slightly from \cite{Enriquez:EllAss} in that we define the elliptic KZB associator using the ``elliptic transport isomorphism'' of Brown--Levin. This definition is analogous to the definition of the Drinfeld associator using parallel transport along the KZ-connection \cite{DG}. We also recall an important result of Enriquez (cf. \cite{Enriquez:EllAss}, \S 6) which describes the variation of the elliptic KZB associator in the modulus of the once-punctured elliptic curve.


\subsection{Definition via the transport function} \label{ssec:4.1}

In Section \ref{ssec:2.3}, we have defined a transport function $T^{\rm KZB}_{\rho_2,\rho_1}$ on a once-punctured elliptic curve for any choice of base points $\rho_1,\rho_2$ (possibly tangential), using the elliptic KZB connection. We now specialize these base points to be $\pm \overrightarrow{v}_0$, where $\overrightarrow{v}_0$ is the tangent vector $-(2\pi i)^{-1}\frac{\partial}{\partial \xi}$ at $0 \in E_{\tau}$. Note that under the isomorphism $E_{\tau} \cong \bC^{\times}/q^\bZ$, we have $\overrightarrow{v}_0=-\frac{\partial}{\partial z}$, where $z=e^{2\pi i\xi}$. In particular, $\overrightarrow{v}_0$ is defined over $\bZ$ on the Tate curve.

Consider now the paths $\alpha,\beta \in \pi_1(E_{\tau}^{\times};-\overrightarrow{v}_0,\overrightarrow{v}_0)$ which are the images of, respectively, the (open) straight-line paths $(0,1)$ and $(0,\tau)$ under the projection $\bC \setminus (\bZ+ \bZ\tau) \rightarrow E_{\tau}^{\times}$, where the path $(0,\tau)$ is additionally composed with a half-circle in the positive direction around $\tau$. Therefore (after ignoring the $-(2\pi i)^{-1}$-prefactor), the paths $\alpha,\beta$ look like in Figure 1 below (cf. \cite{Enriquez:EllAss}, p.550).
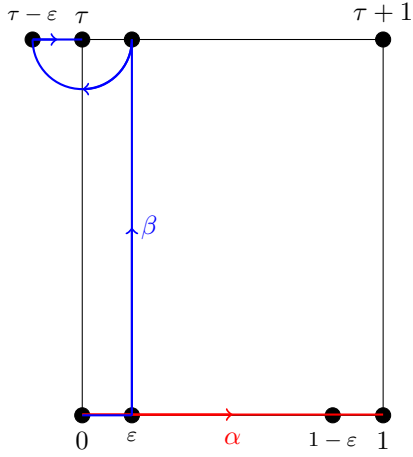
\begin{figure}[h]
	\begin{tikzpicture}
	\draw (0,0) -- (4,0);
	\draw (0,0) -- (0,5);
	\draw (4,0) -- (4,5);
	\draw (0,5) -- (4,5);
	
	\node [below] at (0,-0.1) {$0$};
	\node [below] at (4,-0.1) {$1$};
	\node [above] at (0,5.1) {$\tau$};
	\node [above] at (4,5.1) {$\tau+1$};
	\node [red,below] at (2,-0.1) {$\alpha$};
	\node [blue,right] at (0.66,2.5) {$\beta$};
	\node [below] at (0.66,-0.1) {\footnotesize$\varepsilon$};
	\node [below] at (3.33,-0.1) {\footnotesize$1-\varepsilon$};
	\node [above] at (-0.66,5.1) {\footnotesize$\tau-\varepsilon$};
	
	\draw [fill] (0,0) circle [radius=0.1];
	\draw [fill] (0.66,0) circle [radius=0.1];
	\draw [fill] (3.33,0) circle [radius=0.1];
	\draw [fill] (4,0) circle [radius=0.1];
	\draw [fill] (0,5) circle [radius=0.1];
	\draw [fill] (4,5) circle [radius=0.1];
	\draw [fill] (-0.66,5) circle [radius=0.1];
	\draw [fill] (0.66,5) circle [radius=0.1];
	
	\draw [->,red,thick] (0,0.01) -- (2,0.01);
	\draw [red,thick] (0,0.01) -- (4,0.01);
	\draw [blue,thick] (0,0) -- (0.66,0);
	\draw [->,blue,thick] (0.66,0) -- (0.66,2.5);
	\draw [blue,thick] (0.66,0) -- (0.66,5);
	\draw [->,blue,thick] (-0.66,5) -- (-0.33,5);
	\draw [blue,thick] (-0.66,5) -- (0,5);
	\draw [->,blue,thick] (0.66,5) arc (0:-90:0.66);
	\draw [blue,thick] (-0.66,5) arc (-180:0:0.66);
	
	\end{tikzpicture}
	\caption{The paths $\alpha$ and $\beta$.}
\end{figure}
\begin{dfn}[\cite{Enriquez:EllAss}, \S 6.2]
The \textit{elliptic KZB associator} is the tuple $(A(\tau),B(\tau))$, where
\begin{equation} \label{eqn:AB}
A(\tau):=T^{\rm KZB}_{-\overrightarrow{v}_0,\overrightarrow{v}_0}(\alpha), \quad B(\tau):=T^{\rm KZB}_{-\overrightarrow{v}_0,\overrightarrow{v}_0}(\beta)
\end{equation}
are the images of the paths $\alpha$ and $\beta$ under the transport map $T^{\rm KZB}_{-\overrightarrow{v}_0,\overrightarrow{v}_0}$. 
\end{dfn}
\begin{rmk}
The definition of the elliptic KZB associator given here is not exactly the same as the one given in \cite{Enriquez:EllAss}, but equivalent. Using the elliptic transport map, Enriquez definition is
\begin{equation}
A^{\rm Enr}(\tau):=T^{\rm KZB}_{\overrightarrow{v}_0}(\alpha), \quad B^{\rm Enr}(\tau):=T^{\rm KZB}_{\overrightarrow{v}_0}(\beta).
\end{equation}
Explicitly, the relation between the two versions is given by
\begin{equation}
A(\tau)=e^{-\pi i[\sa,\sb]}A^{\rm Enr}(\tau), \quad B(\tau)=e^{\pi i[\sa,\sb]}B^{\rm Enr}(\tau).
\end{equation}
\end{rmk}

\subsection{Variation in the modulus} \label{ssec:4.2}

An important property of the elliptic KZB associator is that it satisfies a linear differential equation, which relates it to iterated Eisenstein integrals and the special derivations $\varepsilon_{2k}$ reviewed in Section \ref{sec:2}. The boundary condition of this differential equation establishes a relation between the series $A(\tau)$, $B(\tau)$ and the Drinfeld associator $\Phi$. More precisely, we have the following theorem, due to Enriquez.
\begin{thm}[\cite{Enriquez:Emzv}, \S 5.2] \label{thm:Diffeq}
We have
\begin{equation}
A(\tau)=g(\tau)(A_{\infty}), \quad B(\tau)=g(\tau)(B_{\infty}),
\end{equation}
where
\begin{equation}
g(\tau)=\sum (-2\pi i)^n\cG(2k_1,\ldots,2k_n;\tau)\cdot (\varepsilon_{2k_1}\circ \ldots \circ \varepsilon_{2k_n}),
\end{equation}
the sum being over all multi-indices $(k_1,\ldots,k_n) \in \bZ_{\geq 0}^n$, for $n \geq 0$, 
and
\begin{align}
A_{\infty}&=e^{\pi i\iota(\sx_1)}\Phi(\iota(\sx_0)),\iota(\sx_1))e^{2\pi i\iota(\sx_0)}\Phi(\iota(\sx_0),\iota(\sx_1))^{-1},\\
B_{\infty}&=\Phi(\iota(\sx_{\infty})),\iota(\sx_1))e^a\Phi(\iota(\sx_0),\iota(\sx_1))^{-1},
\end{align}
where $\iota: \fp(U) \rightarrow \fp(E_{\tau}^{\times})$ is the morphism of Remark \ref{rmk:elldepth}.
\end{thm}
The element $g(\tau)$ defines an automorphism of $\exp \fp(E_{\tau}^{\times})$.
Letting
\begin{alignat}{3}
&\fA(\tau)&&:=\log(A(\tau)), \quad \fB(\tau)&&:=\log(B(\tau)),\\
&\fA_{\infty}&&:=\log(A_{\infty}), \quad \fB_{\infty}&&:=\log(B_{\infty}),\\
\end{alignat}
we also have
\begin{equation}
\fA(\tau)=g(\tau)(\fA_{\infty}), \quad \fB(\tau)=g(\tau)(\fB_{\infty}),
\end{equation}
since $g(\tau)$ commutes with exponential and logarithm functions.

The next corollary follows immediately from Proposition \ref{prop:metabaction}.
\begin{cor} \label{cor:gmetab}
Let $g(\tau)^{\rm met-ab}$ be the image of $g(\tau)$ in $\End(\fp(E_{\tau}^{\times})_{\bC}^{\rm met-ab})$. We have
\begin{align}
g(\tau)^{\rm met-ab}= &\sum_{n \geq 0}(-2\pi i)^n\cG(\{0\}_n;\tau)\cdot\varepsilon_0^n\\
&+\sum_{n\geq 0, \, k \geq 1}(-2\pi i)^{n+1}\cG(\{0\}_n,2k;\tau)\cdot\Big(\varepsilon_0^n \circ \varepsilon_{2k}\Big)\\
&+\sum_{k,n\geq 1}(-2\pi i)^{n+1}\cG(\{0\}_{n-1},2k,0;\tau)\cdot\Big(\varepsilon_0^{n-1}\circ \varepsilon_{2k}\circ \varepsilon_0\Big).
\end{align}
\end{cor}

\begin{rmk}
The pair $(A_{\infty},B_{\infty})$ is the image of the Drinfeld associator under the natural map (\cite{Enriquez:EllAss}, \S 4.5)
\begin{equation}
\underline{M}(\bC)\rightarrow \underline{Ell}(\bC),
\end{equation}
where $\underline{M}$ is the scheme of classical associators in the sense of \cite{Drinfeld:Gal}, and $\underline{Ell}$ is its elliptic counterpart \cite{Enriquez:EllAss}. A geometric way of interpreting this morphism is via the degeneration of the once-punctured Tate curve to $\bP^1 \setminus \{0,1,\infty\}$ (cf. Remark \ref{rmk:elldepth}).
\end{rmk}

\subsection{Elliptic KZB associator in depth zero} \label{ssec:4.3}

Let $\fA(\tau)^{0}$ be the image of $\fA(\tau)$ in $\gr^0_D\fp(E_{\tau}^{\times})_{\bC}=\fp(E_{\tau}^{\times})/[\fp(E_{\tau}^{\times}),\fp(E_{\tau}^{\times})]$, and likewise let $\fB(\tau)^{(0)}$ be the image of $\fB(\tau)$ in $\gr^0_D\fp(E_{\tau}^{\times})_{\bC}$.

The following proposition shows that $\fA(\tau)^0$ and $\fB(\tau)^0$ precisely retrieve the periods of $H^1(E_{\tau}^{\times})$.
\begin{prop} \label{prop:ab}
We have
\begin{equation}
\fA(\tau)^{(0)}=2\pi i\sb, \quad \fB(\tau)^{(0)}=\sa+2\pi i\tau \sb.
\end{equation}
\end{prop}
\begin{prf}
We only prove the result for $\fA(\tau)^{(0)}$, the formula for $\fB(\tau)^{(0)}$ is proved analogously. By Theorem \ref{thm:Diffeq}, we know that $A(\tau)=g(\tau)(A_{\infty})$, and since $g(\tau)$ is an automorphism, we also have
\begin{equation}
\fA(\tau)=\log(A(\tau))=g(\tau)(\log(A_{\infty}))=g(\tau)(\fA_{\infty}).
\end{equation}
On the other hand, it follows directly from the explicit formula for $A_{\infty}$ given in Theorem \ref{thm:Diffeq} that 
\begin{equation}
\fA_{\infty} \equiv 2\pi i\sb \mod D^1\fp(E_{\tau}^{\times}),
\end{equation}
since $\iota(\sx_0) \equiv \sb \mod D^1\fp(E_{\tau}^{\times})$ and $\iota(\sx_1) \equiv 0 \mod D^1\fp(E_{\tau}^{\times})$. But as every derivation $\varepsilon_{2k}$ annihilates $\sb$, we finally get $\fA(\tau)^{(0)}=g(\tau)(2\pi i\sb)=2\pi i\sb$.
\end{prf}
\begin{rmk} \label{rmk:alt}
Proposition \ref{prop:ab} could have also been proved directly without recourse to Enriquez' Theorem \ref{thm:Diffeq}, using that $\omega_{\rm KZB} \equiv \dd r\cdot \sa+2\pi i\dd\xi \cdot \sb \mod D^1\fp(E_{\tau}^{\times})$.
\end{rmk}

\section{The meta-abelian elliptic KZB associator} \label{sec:5}

In this section, we compute the image of $\fA(\tau)$ and $\fB(\tau)$ in the meta-abelian quotient $\fp(E_{\tau}^{\times})_{\bC}^{\rm met-ab}$ of $\fp(E_{\tau}^{\times})_{\bC}$. The strategy is to use Theorem \ref{thm:Diffeq} which yields that
\begin{equation} \label{eqn:crucial2}
\fA(\tau)=g(\tau)(\fA_{\infty}), \quad \fB(\tau)=g(\tau)(\fB_{\infty})
\end{equation}
and then to compute the images of $\fA_{\infty}$ and $\fB_{\infty}$ in the meta-abelian quotient separately. This is done in Section \ref{ssec:5.1}. In Section \ref{ssec:5.2}, we then compute the action of $g(\tau)$ on the meta-abelian quotient. The two computations are then combined in Section \ref{ssec:5.3} to yield our formula for $\fA(\tau)^{\rm met-ab}$ and $\fB(\tau)^{\rm met-ab}$.

\subsection{The arithmetic piece: periods of Eisenstein series} \label{ssec:5.1}

Let $\fA^{\rm met-ab}_{\infty}$ (resp. $\fB^{\rm met-ab}_{\infty}$) be the image of $\fA_{\infty}$ (resp. the image of $\fB_{\infty}$) in the meta-abelian quotient $\fp(E_{\tau}^{\times})_{\bC}^{\rm met-ab} \cong \gr^0_D\fp(E_{\tau}^{\times})_{\bC} \ltimes \gr^1_D\fp(E_{\tau}^{\times})_{\bC}$, so that we can write
\begin{equation}
\fA_{\infty}^{\rm met-ab}=\fA_{\infty}^{(0)}+\fA_{\infty}^{(1)}, \quad \fB_{\infty}^{\rm met-ab}=\fB_{\infty}^{(0)}+\fB_{\infty}^{(1)}.
\end{equation}
The computation of the depth zero component was already carried out in Proposition \ref{prop:ab} so that it remains to compute the depth one contribution. For this, we need a short lemma about the Drinfeld associator.
\begin{lem} \label{lem:Drinf}
Let $\varphi(\sx_0,\sx_1):=\log(\Phi(\sx_0,\sx_1))$. Then
\begin{equation}
\varphi(\iota(\sx_0),\iota(\sx_1)) \equiv -\sum_{n \geq 2}\zeta(n)\ad^{n-1}(\sb)([\sa,\sb]) \mod D^2\fp(E_{\tau}^{\times})_{\bC},
\end{equation}
where $\iota(\sx_0)=\frac{\ad(\sa)}{e^{\ad(\sa)}-1}(\sb)$ and $\iota(\sx_1)=[\sa,\sb]$ (cf. Remark \ref{rmk:elldepth}). In particular, we have $\varphi(\iota(\sx_0),\iota(\sx_1)) \in D^1\fp(E_{\tau}^{\times})_{\bC}$.
\end{lem}
\begin{prf}
It is well-known (cf. \cite{DG}, \S 6.7) that
\begin{equation}
\varphi(\sx_0,\sx_1) \equiv -\sum_{n=2}^{\infty}\zeta(n)\ad^{n-1}(\sx_0)(\sx_1).
\end{equation}
Applying $\iota$ to both sides, we get the result.
\end{prf}
\begin{thm} \label{thm:arithmetic}
We have
\begin{align}
\fA_{\infty}^{(1)}&=2\pi i\left(c(\sU)-\frac{2\pi i}{4} \sV+\sum_{n \geq 3, {\rm odd}}\zeta(n)\sV^n\right), \label{eqn:constA}\\
\fB_{\infty}^{(1)}&=-2\pi i\left(c(2\pi i\sV)-\sU c(\sU)c(2\pi i\sV)\right)+\sum_{n \geq 3, \,{\rm odd}}\zeta(n)\sU\sV^{n-1},\label{eqn:constB}
\end{align}
where $c(x)=\frac{1}{e^x-1}+\frac 12-\frac 1x=\sum_{k=2}^{\infty}\frac{B_k}{k!}x^{k-1}$.
\end{thm}

\begin{prf}
By Theorem \ref{thm:Diffeq}, we know that
\begin{equation}
\fA_{\infty}=\log(e^{\pi i\iota(\sx_1)}\Phi(\iota(\sx_0),\iota(\sx_1))e^{2\pi i\iota(\sx_0)}\Phi(\iota(\sx_0),\iota(\sx_1))^{-1}).
\end{equation}
Using a ``truncated'' version of the Baker--Campbell--Hausdorff formula (cf. \cite{Reutenauer}, Corollary 3.24) and Lemma \ref{lem:Drinf}, we get
\begin{align}
\mathfrak{S}&:=\log(e^{\pi i\iota(\sx_1)}\Phi(\iota(\sx_0),\iota(\sx_1)))\\
&\equiv \varphi(\iota(\sx_0),\iota(\sx_1))+\sum_{k \geq 0}\frac{B_k}{k!}\ad^k(\varphi(\iota(\sx_0),\iota(\sx_1)))(\pi i\iota(\sx_1))\\
&\equiv \varphi(\iota(\sx_0),\iota(\sx_1))+\pi i\iota(\sx_1) \mod D^2\fp(E_{\tau}^{\times})_{\bC}.\label{eqn:S}
\end{align}
Similarly, since $\iota(\sx_0) \equiv b \mod D^1\fp(E_{\tau}^{\times})_{\bC}$, we get
\begin{align}
\mathfrak{T}&:=\log(e^{2\pi i\iota(\sx_0)}\Phi(\iota(\sx_0),\iota(\sx_1))^{-1})\\
&\equiv -\log(\Phi(\iota(\sx_0),\iota(\sx_1))e^{-2\pi i\iota(\sx_0)})\\
&\equiv 2\pi i\iota(\sx_0)-\sum_{k \geq 0}\frac{B_k}{k!}\ad^k(-\sb)(\varphi(\iota(\sx_0),\iota(\sx_1))) \mod D^2\fp(E_{\tau}^{\times})_{\bC}.\label{eqn:T}
\end{align}
Combining \eqref{eqn:S} and \eqref{eqn:T} and again applying \cite{Reutenauer}, Corollary 3.24, we get
\begin{align}
\fA_{\infty}&\equiv \mathfrak{T}+\sum_{n \geq 0}\frac{B_n}{n!}\ad^n(\mathfrak{T})(\mathfrak{S})\notag\\
&\begin{aligned}\equiv 2\pi i\iota(\sx_0)&-\sum_{k \geq 0}\frac{B_k}{k!}\ad^k(-\sb)(\varphi(\iota(\sx_0),\iota(\sx_1)))+\varphi(\iota(\sx_0),\iota(\sx_1))+\pi i\iota(\sx_1)\\
&+\sum_{n \geq 1}\frac{B_n}{n!}\ad^n(2\pi i\iota(\sx_0))(\varphi(\iota(\sx_0),\iota(\sx_1))+\pi i\iota(\sx_1))
\end{aligned}
\\
&\begin{aligned}\equiv 2\pi i\iota(\sx_0)+\pi i\iota(\sx_1)&-\sum_{k \geq 1}\frac{B_k}{k!}\Big((-1)^k-1 \Big)\ad^k(\sb)(\varphi(\iota(\sx_0),\iota(\sx_1)))\\
&+\sum_{n \geq 1}\frac{B_n}{n!}\ad^n(\sb)(\pi i\iota(\sx_1))\end{aligned}\\
&\begin{aligned}\equiv 2\pi i\sb&+2\pi i\sum_{k \geq 2}\frac{B_k}{k!}\ad^{k-1}(\sa)([\sa,\sb])-\ad(\sb)(\varphi(\iota(\sx_0),\iota(\sx_1)))\\&+
\frac{2\pi i}{2}\sum_{n \geq 1}\frac{B_n(2\pi i)^n}{n!}\ad^n(\sb)([\sa,\sb])\mod D^2\fp(E_{\tau}^{\times})_{\bC},\label{eqn:Aeq1}
\end{aligned}
\end{align}
where in the last line, we have used that $B_1=-\frac 12$ and that $B_{2n+1}=0$ for all $n \geq 1$. Using Lemma \ref{lem:Drinf} together with Euler's formula $-\frac{\zeta(k)}{(-2\pi i)^{k}}=\frac{B_k}{2k!}$ for $k \geq 2$ even, it follows that \eqref{eqn:Aeq1} equals
\begin{equation} \label{eqn:Aeq2}
2\pi i\left(\sb+\sum_{k \geq 2}\frac{B_k}{k!}\ad^{k-1}(\sa)([\sa,\sb])-\frac{2\pi i}{4}\ad(\sb)([\sa,\sb])+\sum_{n \geq 3, {\rm odd}}\zeta(n)\ad^n(\sb)([\sa,\sb])\right).
\end{equation}
Under the substitution $\ad^k(\sa)\ad^l(\sb)([\sa,\sb]) \mapsto \sU^k\sV^l$ (cf. \eqref{eqn:isom2}), \eqref{eqn:constA} now follows immediately from \eqref{eqn:Aeq2} (the $2\pi i\sb$-term belongs to $\fA_{\infty}^{(0)}$ and does not contribute to $\fA_{\infty}^{(1)}$).
The calculation of $\fB_{\infty}^{(1)}$ is very similar, so we will omit some details. First, by definition
\begin{equation}
\fB_{\infty}=\log(\Phi(\iota(\sx_{\infty}),\iota(\sx_1))e^{\sa}\Phi(\iota(\sx_0),\iota(\sx_1))^{-1}),
\end{equation}
where $\sx_{\infty}:=-\sx_0-\sx_1$.
Furthermore,
\begin{align} 
\mathfrak{T}&:=\log(e^\sa\Phi(\iota(\sx_0),\iota(\sx_1))^{-1})\\
&\equiv -\log(\Phi(\iota(\sx_0),\iota(\sx_1))e^{-\sa})\\
&\equiv \sa-\sum_{k \geq 0}\frac{B_k}{k!}\ad^k(-\sa)(\varphi(\iota(\sx_0),\iota(\sx_1))) \mod D^2\fp(E_{\tau}^{\times})_{\bC}.\label{eqn:B2}\\
\end{align}
We obtain
\begin{align}
\fB_{\infty} &\equiv \log(\Phi(\iota(\sx_{\infty}),\iota(\sx_1))e^\mathfrak{T})\\
&\equiv \mathfrak{T}+\sum_{k \geq 0}\frac{B_k}{k!}\ad^k(\sa)(\varphi(\iota(\sx_{\infty}),\iota(\sx_1))) \mod D^2\fp(E_{\tau}^{\times})_{\bC},
\end{align}
where the last equality follows from the fact that $\mathfrak{T} \equiv \sa \mod D^1\fp(E_{\tau}^{\times})_{\bC}$. A short calculation shows that 
\begin{align}
&\mathfrak{T}+\sum_{k \geq 0}\frac{B_k}{k!}\ad^k(\sa)(\varphi(\iota(\sx_{\infty}),\iota(\sx_1)))\\
&\equiv \sa-\sum_{k \geq 0}\frac{B_k}{k!}(-1)^k\ad^k(\sa)(\varphi(\iota(\sx_0),\iota(\sx_1)))+\sum_{k \geq 0}\frac{B_k}{k!}\ad^k(\sa)(\varphi(\iota(\sx_{\infty}),\iota(\sx_1)))\\
& \equiv \sa+\sum_{k\geq 0}\frac{B_k}{k!}\ad^k(\sa)\Big( \varphi(\iota(\sx_{\infty}),\iota(\sx_1))-(-1)^k\varphi(\iota(\sx_0),\iota(\sx_1)) \Big) \mod D^2\fp(E_{\tau}^{\times})_{\bC}\\
\label{eqn:B3}
\end{align}
The term in brackets is equal to
\begin{equation}
\begin{cases}
\displaystyle 2\sum_{n\geq 2, \, \rm{even}}\zeta(n)\ad^{n-1}(\sb)([\sa,\sb]) &\mbox{if $k$ is even}\\
\displaystyle-2\sum_{n\geq 3 \, \rm{odd}}\zeta(n)\ad^{n-1}(\sb)([\sa,\sb]) & \mbox{if $k$ is odd.}
\end{cases}
\end{equation}
Again using that $\zeta(k)=-\frac{B_k(2\pi i)^k}{2k!}$, if $k \geq 2$ is even, we obtain that \eqref{eqn:B3} equals
\begin{align}
\sa-\sum_{n \geq 2}\frac{B_n(2\pi i)^n}{n!}\ad^{n-1}(\sb)([\sa,\sb])&-\sum_{k,n \geq 2}\frac{B_kB_n(2\pi i)^n}{k!n!}\ad^k(\sa)\ad^{n-1}(\sb)([\sa,\sb])\\
&+\begin{aligned}\sum_{n\geq 3, \, \rm{odd}}\zeta(n)\ad(\sa)&\ad^{n-1}(\sb)([\sa,\sb]) \\
& \mod D^2\fp(E_{\tau}^{\times})_{\bC}.
\end{aligned}\label{eqn:B4}
\end{align}
The first term $\sa$ belongs to $\fB_{\infty}^{(0)}$, and does not contribute to $\fB^{(1)}_{\infty}$. Applying the isomorphism \eqref{eqn:isom2} to the remaining terms in \eqref{eqn:B4}, we obtain the desired result \eqref{eqn:constB}.
\end{prf}
The series $\fA_{\infty}^{(1)}$ and $\fB_{\infty}^{(1)}$ are closely related to the extended period polynomials of Eisenstein series $r_{G_{2k}}(\sX,\sY)$ \cite{Zagier:Periods}. Precisely, for $k \geq 2$, one has
\begin{equation} \label{eqn:periodpolynomial}
r_{G_{2k}}(\sX,\sY)=\omega_{G_{2k}}^+P_{G_{2k}}(\sX,\sY)^++\omega_{G_{2k}}^-P_{G_{2k}}(\sX,\sY)^-,
\end{equation}
where
\begin{align}
P_{G_{2k}}(\sX,\sY)^+&=\sX^{2k-2}-\sY^{2k-2}\\
P_{G_{2k}}(\sX,\sY)^-&=\sum_{-1\leq n \leq 2k-1}\frac{B_{n+1}B_{2k-n-1}}{(n+1)!(2k-1-n)!}\sX^n\sY^{2k-2-n}
\end{align}
and $\omega_{G_{2k}}^-=-\frac{(2k-2)!}{2}$, $\omega_{G_{2k}}^+=\frac{\zeta(2k-1)}{(2\pi i)^{2k-1}}\omega_{G_{2k}}^-$ (the ``periods'' of $G_{2k}$). Now let
\begin{equation}
\widetilde{\fA}(\sU,\sV)=\frac{1}{\sV}\fA_{\infty}^{(1)}(\sU,\sV), \quad \widetilde{\fB}(\sU,\sV)=\frac{1}{\sU}\fB_{\infty}^{(1)}(\sU,\sV).
\end{equation}
These are formal Laurent series in the variables $\sU$ and $\sV$. In general, if $f(\sU,\sV)$ is a formal Laurent series, we denote by $f(\sU,\sV)_k$ its homogeneous component of degree $k$ and $f(\sU,\sV)^\pm:=\frac{f(\sU,\sV)\pm f(-\sU,\sV)}{2}$. Comparing now \eqref{eqn:periodpolynomial} with Theorem \ref{thm:arithmetic}, we get
\begin{cor} \label{cor:arithmetic}
We have
\begin{align}
r_{G_{2k}}(\sU,\sV)=
\frac{\omega_{G_{2k}}^-}{2\pi i}\Bigg[\widetilde{\fA}(\overline{\sU},\sV)_{2k-2}^++\widetilde{\fB}(\sV,\overline{\sU})_{2k-2}^+
-\widetilde{\fA}(\overline{\sU},\sV)_{2k-2}^--\widetilde{\fB}(\overline{\sU},\sV)_{2k-2}^-\Bigg],
\end{align}
where $\overline{\sU}=\frac{\sU}{2\pi i}$.
\end{cor}

\subsection{The geometric piece: special values of elliptic polylogarithms} \label{ssec:5.2}

Recall from Section \ref{ssec:4.2} the definition of the automorphism $g(\tau): \exp \fp(E_{\tau}^{\times})_{\bC} \rightarrow \exp \fp(E_{\tau}^{\times})_{\bC}$. It naturally extends to the topological enveloping algebra $\bQ\langle\!\langle \sa,\sb\rangle\!\rangle$ of $\fp(E_{\tau}^{\times})_{\bC}$.

In this section, we compute the images of $g(\tau)(\sa)$, $g(\tau)(\sb)$ in the meta-abelian quotient $\fp(E_{\tau}^{\times})_{\bC}^{\rm met-ab}$ of $\fp(E_{\tau}^{\times})_{\bC}$, and relate the result to special values of Beilinson--Levin's elliptic polylogarithms \cite{BeiLev,Levin:Compositio}.
\begin{thm} \label{thm:geometric}
Let $\sW=\frac{\sU}{2\pi i}+\tau \sV$. We have
\begin{equation} \label{eqn:geom1}
g(\tau)(\sa)^{\rm met-ab}=\sa+2\pi i\tau \sb-2\pi i\sW\sum_{k=1}^{\infty}\frac{2}{(2k-2)!} \int_{\tau}^{\overrightarrow{1}_{\infty}}\underline{G}_{2k},
\end{equation}
and
\begin{equation} \label{eqn:geom2}
g(\tau)(\sb)^{\rm met-ab}=2\pi i\sb-2\pi i\sV\sum_{k=1}^{\infty}\frac{2}{(2k-2)!} \int_{\tau}^{\overrightarrow{1}_{\infty}}\underline{G}_{2k},
\end{equation}
where $\underline{G}_{2k}=(2\pi i)^{2k-1}G_{2k}(z)(\sW-z\sV)^{2k-2}\dd z$.
\end{thm}

\begin{prf}
By Corollary \ref{cor:gmetab}, we have
\begin{align}
g(\tau)(\sa)^{\rm met-ab}& =\sa+2\pi i\tau \sb+\sum_{n\geq 0, \, k\geq 1}(-2\pi i)^{n+1}\cG(\{0\}_n,2k;\tau)\Big(\varepsilon_0^n \circ \varepsilon_{2k}\Big)(\sa)\\
&+\sum_{k,n \geq 1}(-2\pi i)^{n+1}\cG(\{0\}_{n-1},2k,0;\tau)\Big(\varepsilon_0^{n-1} \circ \varepsilon_{2k} \circ \varepsilon_0\Big)(\sa)\\
&= \sa+2\pi i\tau \sb+\sum_{n\geq 0, \, k\geq 1}\frac{2(-2\pi i)^{n+1}}{(2k-2)!}\cG(\{0\}_n,2k;\tau)\varepsilon_0^n(\ad^{2k-1}(\sa)([\sa,\sb]))\\
&-\sum_{k,n \geq 1}\frac{2(-2\pi i)^{n+1}}{(2k-2)!}\cG(\{0\}_{n-1},2k,0;\tau)\varepsilon_0^{n-1}(\ad^{2k-2}(\sa)\ad(\sb)([\sa,\sb])), \label{eqn:geom3}
\end{align}
Using the isomorphism of Proposition \ref{prop:metablie} together with Proposition \ref{prop:metabaction} and Proposition \ref{prop:IterEisEichler}, we see that \eqref{eqn:geom3} equals
\begin{align}
&\sa+2\pi i\tau\sb-\sum_{n\geq 0, \, k\geq 1}\frac{2(2\pi i)^{n+1}}{(2k-2)!n!}I_n(G_{2k};\tau)\left(\sV\frac{\partial}{\partial \sU} \right)^n\sU^{2k-1}\\
&-\sum_{k,n \geq 1}\frac{2(2\pi i)^{n+1}}{(2k-2)!(n-1)!}\Big(\tau I_{n-1}(G_{2k};\tau)-I_n(G_{2k};\tau) \Big) \left(\sV\frac{\partial}{\partial \sU} \right)^{n-1} \sU^{2k-2}\sV.
\end{align}
Now we apply the differential operator $\sV\frac{\partial }{\partial \sU}$ and split the first and the last sum to obtain
\begin{align}
g(\tau)(\sa)^{\rm met-ab}&=\sa+2\pi i\tau \sb-\sum_{k\geq 1}\frac{2(2\pi i)}{(2k-2)!}I_0(G_{2k};\tau)\sU^{2k-1}\\
&-\sum_{k,n\geq 1}\frac{2(2k-1)(2\pi i)^{n+1}}{(2k-1-n)!n!}I_n(G_{2k};\tau)\sU^{2k-1-n}\sV^n\\
&-2\pi i\tau\sum_{k,n \geq 1}\frac{2(2\pi i)^n}{(2k-1-n)!(n-1)!} I_{n-1}(G_{2k};\tau)\sU^{2k-1-n}\sV^{n-1}\\
&+\sum_{k,n \geq 1}\frac{2(2\pi i)^{n+1}}{(2k-1-n)!(n-1)!}I_n(G_{2k};\tau)\sU^{2k-1-n}\sV^{n-1}.
\end{align}
From the definition of $I_n(G_{2k};\tau)$, it is easy to see that the third sum equals
\begin{equation}
-2\pi i\tau \sV\sum_{k=1}^{\infty}\frac{2(2\pi i)}{(2k-2)!}\int_{\tau}^{\overrightarrow{1}_{\infty}}G_{2k}(z)\Big(\sU+2\pi i(\tau-z)\sV\Big)^{2k-2}\dd z.
\end{equation}
On the other hand, the first, second and fourth sum give
\begin{equation}
-\sU\sum_{k=1}^{\infty}\frac{2(2\pi i)}{(2k-2)!}\int_{\tau}^{\overrightarrow{1}_{\infty}}G_{2k}(z)\Big(\sU+2\pi i(\tau-z)\sV\Big)^{2k-2}\dd z.
\end{equation}
Combining the two equations and setting $\sW=\frac{\sU}{2\pi i}+\tau \sV$, the first equality \eqref{eqn:geom1} follows. Since $g(\tau)$ is uniquely determined by its value on $e^\sa$, the second statement \eqref{eqn:geom2} follows from the first, but can also be proved directly along similar lines.
\end{prf}
We now give the relation to special values of elliptic polylogarithms. Following the notation of \cite{Levin:Compositio}, we let $\Xi(\xi,\tau;\sX,\sY)$ be the (modified) generating series of elliptic polylogarithms $\Lambda_{m,n}(\xi,\tau)$. These are holomorphic functions on the universal covering of the once-punctured elliptic curve $E_{\tau}^{\times}$, which are obtained by averaging the (Debye) polylogarithms along the spiral $q^{\bZ}$. Let
\begin{equation}
\Xi^*(0,\tau;\sX,\sY):=(\Xi(\xi,\tau;\sX,\sY)-\frac{1}{2\pi i}\log(2\pi i\xi))\vert_{\xi=0}
\end{equation}
be its (regularized) special value at the zero section of the elliptic curve. It has been shown in \cite{Levin:Compositio}, Theorem 4.1 that
\begin{equation} \label{eqn:Levin}
\Xi^*(0,\tau;\sX,\sY)=\frac{-\tau}{X(X-\tau Y)}+\sum_{k=2}^{\infty}(-1)^{k-1}(k-1)\cE_k,
\end{equation}
where for $k\geq 2$, $\cE_k$ is the indefinite integral of $E_k(\tau)(X-\tau Y)^{k-2}\dd\tau$ with $E_k(\tau)=\frac{2(2\pi i)^k}{(k-1)!}G_k(\tau)=\sum_{(m,n) \in \bZ^2 \setminus \{(0,0)\}}\frac{1}{(m\tau+n)^k}$ the classical Eisenstein series of weight $k$. The constants of integration in the indefinite integrals can be retrieved uniformly as the (regularized) special value of $\Xi^*(0,\tau;X,Y)$ at $\tau=i\infty$, which is straightforwardly computed from the definitions and is given explicitly by
\begin{align} \label{eqn:constellpol}
\Xi^*(0,i\infty;X,Y)=-\sum_{n\geq 2}\frac{\zeta(n)}{(2\pi i)^n}Y^{n-1}+\frac{1}{e^X-1}\left( \frac{1}{e^Y-1}-\frac 1Y \right).
\end{align}
Now comparing \eqref{eqn:Levin} with Theorem \ref{thm:geometric}, we obtain
\begin{cor} \label{cor:geometric}
Let $g(\tau)(a)^{\rm met-ab}-a$, and replace $2\pi ib$ by $(\sW-\tau \sV)^{-1}$. Then
\begin{align}
\frac{g(\tau)(\sa)^{\rm met-ab}-\sa}{-(2\pi i)^2\sW}=\Xi^*(0,\tau;2\pi i\sW,2\pi i\sV)-\Xi^*(0,i\infty;2\pi i\sW,2\pi i\sV),
\end{align}
where $\Xi^*(0,i\infty;X,Y)$ is given in \eqref{eqn:constellpol} above.
\end{cor}

\subsection{Putting the pieces together} \label{ssec:5.3}

We can now complete the computation of $\fA(\tau)^{\rm met-ab}$ and $\fB(\tau)^{\rm met-ab}$ by combining the results of the previous sections.
\begin{thm} \label{thm:arithgeo}
We have
\begin{align}
\fA(\tau)^{\rm met-ab}&=2\pi i\sb+\exp\left(\tau\frac{\partial}{\partial \overline{\sU}}\sV \right)\fA^{(1)}_{\infty}-2\pi i\sV\sum_{k=1}^{\infty}\frac{2}{(2k-2)!}\int_{\tau}^{\overrightarrow{1}_{\infty}}\underline{G}_{2k},
\end{align}
and
\begin{align}
\fB(\tau)^{\rm met-ab}&=\sa+2\pi i\tau \sb+\exp\left(\tau\frac{\partial}{\partial \overline{\sU}}\sV \right)\fB^{(1)}_{\infty}-2\pi i\sW\sum_{k=1}^{\infty}\frac{2}{(2k-2)!}\int_{\tau}^{\overrightarrow{1}_{\infty}}\underline{G}_{2k}.
\end{align}
where $\overline{\sU}=\frac{\sU}{2\pi i}$, $\sW=\overline{\sU}+\tau \sV$ and $\fA^{(1)}_{\infty}$ and $\fB^{(1)}_{\infty}$ are as given in Theorem \ref{thm:arithmetic}
\end{thm}
\begin{prf}
We only prove the first equality, the second one is shown analogously. By Theorem \ref{thm:Diffeq}, we have $\fA(\tau)=g(\tau)(\fA_{\infty})$, hence
\begin{equation}
\fA(\tau)^{\rm met-ab} \equiv g(\tau)(\fA_{\infty}) \mod D^2\fp(E_{\tau}^{\times})_{\bC},
\end{equation}
and from Proposition \ref{prop:metabaction}, we get
\begin{equation}
\fA(\tau)^{\rm met-ab}=g(\tau)(\fA^{(1)}_{\infty})+2\pi ig(\tau)(\sb)^{\rm met-ab}.
\end{equation}
The only derivation which acts non-trivially on $\gr^1_D\fp(E_{\tau}^{\times})_{\bC}$ is $\varepsilon_0$ which itself acts as $-\frac{\partial}{\partial \sU}\sV=\frac{1}{2\pi i}\frac{\partial}{\partial \overline{\sU}}\sV$. Combining this with Theorem \ref{thm:geometric}, we get the result:
\begin{align}
\fA(\tau)^{\rm met-ab}=2\pi i\sb+\exp\left(\tau\frac{\partial}{\partial \overline{\sU}}\sV\right)\fA^{(1)}_{\infty}-2\pi i\sV\sum_{k=1}^{\infty}\frac{2}{(2k-2)!}\int_{\tau}^{\overrightarrow{1}_{\infty}}\underline{G}_{2k}.
\end{align}
\end{prf}

\begin{rmk}
The value for $\fA(\tau)^{\rm met-ab}$ given in Theorem \ref{thm:arithgeo} can be further simplified. To this end, recall from Theorem \ref{thm:arithmetic} that
\begin{equation}
\fA_{\infty}^{(1)}=2\pi i\left(\sum_{k=1}^{\infty}\frac{B_{2k}}{(2k)!}\sU^{2k-1}-\frac{2\pi i}{4} \sV+\sum_{n=3, \, \rm{odd}}\zeta(n)\sV^n\right).
\end{equation}
Therefore
\begin{align}
\exp\left( \tau\frac{\partial}{\partial \overline{\sU}}\sV \right)\fA_{\infty}^{(1)}&=\fA_{\infty}^{(1)}+2\pi i\sum_{k,n\geq 1}\frac{\tau^n}{n!}\frac{B_{2k}}{(2k)!}\left(\frac{\partial}{\partial \overline{\sU}}\sV \right)^n\sU^{2k-1}\\
&=\fA_{\infty}^{(1)}+2\pi i\sV\sum_{k,n \geq 1}\frac{2(2\pi i)^{2k-1}}{(2k-2)!}\Bigg[\frac{\tau^n}{n!}\frac{B_{2k}}{4k}\left(\frac{\partial}{\partial \overline{\sU}}\sV \right)^{n-1}\overline{\sU}^{2k-2}\Bigg]\\
&=\fA_{\infty}^{(1)}+2\pi iV\sum_{k,n \geq 1}\frac{2(2\pi i)^{2k-1}}{(2k-1-n)!}\Bigg[\frac{\tau^n}{n!}\frac{B_{2k}}{4k}\overline{\sU}^{2k-1-n}\sV^{n-1}\Bigg]\\
&=\fA_{\infty}^{(1)}+2\pi i\sV\sum_{k=1}^{\infty}\frac{2(2\pi i)^{2k-1}}{(2k-2)!}\frac{B_{2k}}{4k}\int_0^{\tau}(\overline{\sU}+(\tau-z)\sV)^{2k-2}\dd z.
\end{align}
Note that $-\frac{B_{2k}}{4k}=a_0(G_{2k})$, the zeroth Fourier coefficient of $G_{2k}$. Consequently, we obtain
\begin{align}
\fA(\tau)^{\rm met-ab}&=2\pi i\sb+\fA_{\infty}^{(1)}-2\pi i\sV\sum_{k=1}^{\infty}\frac{2}{(2k-2)!}\int_{\tau}^{i\infty}\underline{G}^0_{2k},
\end{align}
where $\underline{G}^0_{2k}=\underline{G}_{2k}-a_0(\underline{G}_{2k})=\underline{G}_{2k}-(2\pi i)^{2k-1}a_0(G_{2k})(\sW-z\sV)^{2k-2}$, since
\begin{align}
\int_{\tau}^{\overrightarrow{1}_{\infty}}\underline{G}_{2k}=\int_{\tau}^{i\infty}\underline{G}^0_{2k}-\int_0^{\tau}a_0(\underline{G}_{2k}).
\end{align}
\end{rmk}

\begin{bibtex}[\jobname]

@article {And,
    AUTHOR = {Anderson, Greg W.},
     TITLE = {The hyperadelic gamma function},
   JOURNAL = {Invent. Math.},
  FJOURNAL = {Inventiones Mathematicae},
    VOLUME = {95},
      YEAR = {1989},
    NUMBER = {1},
     PAGES = {63--131},
      ISSN = {0020-9910},
     CODEN = {INVMBH},
   MRCLASS = {11S80 (11G20)},
  MRNUMBER = {969414},
MRREVIEWER = {Gerd Faltings},
       DOI = {10.1007/BF01394145},
       URL = {http://dx.doi.org/10.1007/BF01394145},
}

@book {Andre:Motifs,
    AUTHOR = {Andr{\'e}, Yves},
     TITLE = {Une introduction aux motifs (motifs purs, motifs mixtes,
              p\'eriodes)},
    SERIES = {Panoramas et Synth\`eses [Panoramas and Syntheses]},
    VOLUME = {17},
 PUBLISHER = {Soci\'et\'e Math\'ematique de France, Paris},
      YEAR = {2004},
     PAGES = {xii+261},
      ISBN = {2-85629-164-3},
   MRCLASS = {14F42 (11J91 14C25 19E15)},
  MRNUMBER = {2115000 (2005k:14041)},
MRREVIEWER = {Luca Barbieri Viale},
}

@article {Apery:Irrationalite,
	AUTHOR = {Ap{\'e}ry, Roger},
	 TITLE = {Irrationalit{\'e} de {$\zeta(2)$} et {$\zeta(3)$}},
   JOURNAL = {Ast{\'e}risque},
    VOLUME = {61},
      YEAR = {1979},
     PAGES = {11--13},
}
	 
@book {Bailey:Hypergeometric,
    AUTHOR = {Bailey, W. N.},
     TITLE = {Generalized hypergeometric series},
    SERIES = {Cambridge Tracts in Mathematics and Mathematical Physics, No.
              32},
 PUBLISHER = {Stechert-Hafner, Inc., New York},
      YEAR = {1964},
     PAGES = {v+108},
   MRCLASS = {33.20 (40.00)},
  MRNUMBER = {0185155 (32 \#2625)},
}

@article {BR:Irrationalite,
    AUTHOR = {Ball, Keith and Rivoal, Tanguy},
     TITLE = {Irrationalit\'e d'une infinit\'e de valeurs de la fonction
              z\^eta aux entiers impairs},
   JOURNAL = {Invent. Math.},
  FJOURNAL = {Inventiones Mathematicae},
    VOLUME = {146},
      YEAR = {2001},
    NUMBER = {1},
     PAGES = {193--207},
      ISSN = {0020-9910},
   MRCLASS = {11J72 (11M06)},
  MRNUMBER = {1859021},
MRREVIEWER = {F. Beukers},
       DOI = {10.1007/s002220100168},
       URL = {http://dx.doi.org/10.1007/s002220100168},
}

@article {BKT:EllipticPol,
    AUTHOR = {Bannai, Kenichi and Kobayashi, Shinichi and Tsuji, Takeshi},
     TITLE = {On the de {R}ham and {$p$}-adic realizations of the elliptic
              polylogarithm for {CM} elliptic curves},
   JOURNAL = {Ann. Sci. \'Ec. Norm. Sup\'er. (4)},
  FJOURNAL = {Annales Scientifiques de l'\'Ecole Normale Sup\'erieure.
              Quatri\`eme S\'erie},
    VOLUME = {43},
      YEAR = {2010},
    NUMBER = {2},
     PAGES = {185--234},
      ISSN = {0012-9593},
   MRCLASS = {11G55 (11G15 14F30 14G10)},
  MRNUMBER = {2662664 (2011g:11125)},
MRREVIEWER = {Jan Nekov{\'a}{\v{r}}},
}

@ARTICLE{BS:fundamentalLie,
   author = {{Baumard}, S. and {Schneps}, L.},
    title = "{Relations dans l'alg\`ebre de Lie fondamentale des motifs elliptiques mixtes}",
  journal = {ArXiv e-prints},
   volume = {math.NT/1310.5833},
 primaryClass = "math.AG",
 keywords = {Mathematics - Algebraic Geometry},
     year = 2013,
   adsurl = {http://adsabs.harvard.edu/abs/2013arXiv1310.5833B},
  adsnote = {Provided by the SAO/NASA Astrophysics Data System}
}

@incollection {BeiLev,
	AUTHOR = {Be{\u\i}linson, A. and Levin, A.},
	TITLE = {The elliptic polylogarithm},
	BOOKTITLE = {Motives ({S}eattle, {WA}, 1991)},
	SERIES = {Proc. Sympos. Pure Math.},
	VOLUME = {55},
	PAGES = {123--190},
	PUBLISHER = {Amer. Math. Soc., Providence, RI},
	YEAR = {1994},
	MRCLASS = {11G05 (11G09 11G40 14H52 19F27)},
	MRNUMBER = {1265553},
	MRREVIEWER = {J. Browkin},
}

@book {Bloch:Regulators,
    AUTHOR = {Bloch, Spencer J.},
     TITLE = {Higher regulators, algebraic {$K$}-theory, and zeta functions
              of elliptic curves},
    SERIES = {CRM Monograph Series},
    VOLUME = {11},
 PUBLISHER = {American Mathematical Society, Providence, RI},
      YEAR = {2000},
     PAGES = {x+97},
      ISBN = {0-8218-2114-8},
   MRCLASS = {11G55 (11G40 11R70 14G10 19F27)},
  MRNUMBER = {1760901 (2001i:11082)},
MRREVIEWER = {Jan Nekov{\'a}{\v{r}}},
}

@article {MZVmine,
    AUTHOR = {Bl{\"u}mlein, J. and Broadhurst, D. J. and Vermaseren, J. A.
              M.},
     TITLE = {The multiple zeta value data mine},
   JOURNAL = {Comput. Phys. Comm.},
  FJOURNAL = {Computer Physics Communications. An International Journal and
              Program Library for Computational Physics and Physical
              Chemistry},
    VOLUME = {181},
      YEAR = {2010},
    NUMBER = {3},
     PAGES = {582--625},
      ISSN = {0010-4655},
     CODEN = {CPHCBZ},
   MRCLASS = {11M32 (11Y70)},
  MRNUMBER = {2578167 (2011a:11163)},
MRREVIEWER = {Zhonghua Li},
       DOI = {10.1016/j.cpc.2009.11.007},
       URL = {http://dx.doi.org/10.1016/j.cpc.2009.11.007},
}

@book {Bou,
    AUTHOR = {Bourbaki, N.},
     TITLE = {\'{E}l\'ements de math\'ematique. {F}asc. {XXXIV}. {G}roupes
              et alg\`ebres de {L}ie. {C}hapitre {IV}: {G}roupes de
              {C}oxeter et syst\`emes de {T}its. {C}hapitre {V}: {G}roupes
              engendr\'es par des r\'eflexions. {C}hapitre {VI}: syst\`emes
              de racines},
    SERIES = {Actualit\'es Scientifiques et Industrielles, No. 1337},
 PUBLISHER = {Hermann, Paris},
      YEAR = {1968},
     PAGES = {288 pp. (loose errata)},
   MRCLASS = {22.50 (17.00)},
  MRNUMBER = {0240238 (39 \#1590)},
MRREVIEWER = {G. B. Seligman},
}

@article {BK,
    AUTHOR = {Broadhurst, D. J. and Kreimer, D.},
     TITLE = {Association of multiple zeta values with positive knots via
              {F}eynman diagrams up to {$9$} loops},
   JOURNAL = {Phys. Lett. B},
  FJOURNAL = {Physics Letters. B},
    VOLUME = {393},
      YEAR = {1997},
    NUMBER = {3-4},
     PAGES = {403--412},
      ISSN = {0370-2693},
     CODEN = {PYLBAJ},
   MRCLASS = {11M41 (11Z05 57M25 81T18)},
  MRNUMBER = {1435933 (98g:11101)},
MRREVIEWER = {Louis H. Kauffman},
       DOI = {10.1016/S0370-2693(96)01623-1},
       URL = {http://dx.doi.org/10.1016/S0370-2693(96)01623-1},
}

@article {BMS,
	AUTHOR = {Broedel, Johannes and Matthes, Nils and Schlotterer, Oliver},
	TITLE = {Relations between elliptic multiple zeta values and a special
		derivation algebra},
	JOURNAL = {J. Phys. A},
	FJOURNAL = {Journal of Physics. A. Mathematical and Theoretical},
	VOLUME = {49},
	YEAR = {2016},
	NUMBER = {15},
	PAGES = {155--203},
	ISSN = {1751-8113},
	MRCLASS = {11M32 (33E05)},
	MRNUMBER = {3479125},
	DOI = {10.1088/1751-8113/49/15/155203},
	URL = {http://dx.doi.org/10.1088/1751-8113/49/15/155203},
}

@article {BMMS,
	AUTHOR = {Broedel, Johannes and Mafra, Carlos R. and Matthes, Nils and
		Schlotterer, Oliver},
	TITLE = {Elliptic multiple zeta values and one-loop superstring
		amplitudes},
	JOURNAL = {J. High Energy Phys.},
	FJOURNAL = {Journal of High Energy Physics},
	YEAR = {2015},
	NUMBER = {7},
	PAGES = {112, front matter+41},
	ISSN = {1126-6708},
	MRCLASS = {83C47 (83E30)},
	MRNUMBER = {3383100},
	MRREVIEWER = {Farhang Loran},
}

@article {BSS,
    AUTHOR = {Broedel, Johannes and Schlotterer, Oliver and Stieberger,
              Stephan},
     TITLE = {Polylogarithms, multiple zeta values and superstring
              amplitudes},
   JOURNAL = {Fortschr. Phys.},
  FJOURNAL = {Fortschritte der Physik. Progress of Physics},
    VOLUME = {61},
      YEAR = {2013},
    NUMBER = {9},
     PAGES = {812--870},
      ISSN = {0015-8208},
   MRCLASS = {81T30 (11M32 33B30)},
  MRNUMBER = {3104459},
MRREVIEWER = {Giuseppe Nardelli},
       DOI = {10.1002/prop.201300019},
       URL = {http://dx.doi.org/10.1002/prop.201300019},
}

@article{BSST,
      author         = "Broedel, Johannes and Schlotterer, Oliver and Stieberger,
                        Stephan and Terasoma, Tomohide",
      title          = "{All order $\alpha^{\prime}$-expansion of superstring
                        trees from the Drinfeld associator}",
      journal        = "Phys. Rev.",
      volume         = "D89",
      year           = "2014",
      number         = "6",
      pages          = "066014",
      doi            = "10.1103/PhysRevD.89.066014",
      eprint         = "1304.7304",
      archivePrefix  = "arXiv",
      primaryClass   = "hep-th",
      reportNumber   = "DAMTP-2013-23, AEI-2013-195, MPP-2013-120",
      SLACcitation   = "
}

@incollection {Brown:Colombia,
	AUTHOR = {Brown, Francis},
	TITLE = {Iterated integrals in quantum field theory},
	BOOKTITLE = {Geometric and topological methods for quantum field theory},
	PAGES = {188--240},
	PUBLISHER = {Cambridge Univ. Press, Cambridge},
	YEAR = {2013},
	MRCLASS = {81S40 (81T18 81T40)},
	MRNUMBER = {3098088},
	MRREVIEWER = {Roberto Quezada},
}

@article {Brown:MTM,
    AUTHOR = {Brown, Francis},
     TITLE = {Mixed {T}ate motives over {$\mathbb Z$}},
   JOURNAL = {Ann. of Math. (2)},
  FJOURNAL = {Annals of Mathematics. Second Series},
    VOLUME = {175},
      YEAR = {2012},
    NUMBER = {2},
     PAGES = {949--976},
      ISSN = {0003-486X},
   MRCLASS = {11S20 (11M32 14F42)},
  MRNUMBER = {2993755},
MRREVIEWER = {Pierre A. Lochak},
       DOI = {10.4007/annals.2012.175.2.10},
       URL = {http://dx.doi.org/10.4007/annals.2012.175.2.10},
}

@unpublished{Brown:MMV,
     author         = {Brown, Francis},
     title          = {Multiple modular values and the relative completion of the fundamental group of $\mathcal{M}_{1,1}$},
     year           = {2016},
     note			= {arXiv:1407.5167v3},
}

@unpublished{Brown:ICM,
	author         = {Brown, Francis},
	title          = {Motivic periods and $\mathbb{P}^1$ minus three points},
	year           = {2014},
	note			= {proceedings of the ICM (2014)},
}

@article {Brown:depth3,
	AUTHOR = {Brown, Francis},
	TITLE = {ZETA ELEMENTS IN DEPTH 3 AND THE FUNDAMENTAL
		{L}IE ALGEBRA OF THE INFINITESIMAL {T}ATE CURVE},
	JOURNAL = {Forum Math. Sigma},
	FJOURNAL = {Forum of Mathematics. Sigma},
	VOLUME = {5, e1},
	YEAR = {2017},
	PAGES = {56pp},
	ISSN = {2050-5094},
	MRCLASS = {11M32 (11F67)},
	MRNUMBER = {3593496},
	DOI = {10.1017/fms.2016.29},
	URL = {http://dx.doi.org/10.1017/fms.2016.29},
}

@unpublished{Brown:depth,
	author         = {Brown, Francis},
	title          = {Depth-graded motivic multiple zeta values},
	year           = {2013},
	note			= {arXiv:1301.3053},
}

@incollection {Brown:Decomposition,
	AUTHOR = {Brown, Francis},
	TITLE = {On the decomposition of motivic multiple zeta values},
	BOOKTITLE = {Galois-{T}eichm\"uller theory and arithmetic geometry},
	SERIES = {Adv. Stud. Pure Math.},
	VOLUME = {63},
	PAGES = {31--58},
	PUBLISHER = {Math. Soc. Japan, Tokyo},
	YEAR = {2012},
	MRCLASS = {11M32 (13B05 16T15)},
	MRNUMBER = {3051238},
	MRREVIEWER = {Antanas Laurin{\v{c}}ikas},
}

@unpublished{BL:MEP,
     author         = {Brown, Francis. and Levin, Andrey},
     title          = {Multiple elliptic polylogarithms},
	 note   		= {arXiv:1110.6917},
     year           = {2011},
}

@incollection {Car,
    AUTHOR = {Cartier, Pierre},
     TITLE = {A primer of {H}opf algebras},
 BOOKTITLE = {Frontiers in number theory, physics, and geometry. {II}},
     PAGES = {537--615},
 PUBLISHER = {Springer, Berlin},
      YEAR = {2007},
   MRCLASS = {16W30 (01A60 05E05)},
  MRNUMBER = {2290769 (2008b:16059)},
MRREVIEWER = {Ralf Holtkamp},
       DOI = {10.1007/978-3-540-30308-4_12},
       URL = {http://dx.doi.org/10.1007/978-3-540-30308-4_12},
}

@incollection {CEE:KZB,
             AUTHOR = {Calaque, Damien and Enriquez, Benjamin and Etingof, Pavel},
              TITLE = {Universal {KZB} equations: the elliptic case},
          BOOKTITLE = {Algebra, arithmetic, and geometry: in honor of {Y}u. {I}.
                       {M}anin. {V}ol. {I}},
             SERIES = {Progr. Math.},
             VOLUME = {269},
              PAGES = {165--266},
          PUBLISHER = {Birkh\"auser Boston, Inc., Boston, MA},
               YEAR = {2009},
            MRCLASS = {32G34 (11F55 17B37 20C08 32C38)},
           MRNUMBER = {2641173 (2011k:32018)},
         MRREVIEWER = {Gwyn Bellamy},
                DOI = {10.1007/978-0-8176-4745-2_5},
                URL = {http://dx.doi.org/10.1007/978-0-8176-4745-2\_5},
}

@article {Chen:PathIntegrals,
    AUTHOR = {Chen, Kuo Tsai},
     TITLE = {Iterated path integrals},
   JOURNAL = {Bull. Amer. Math. Soc.},
  FJOURNAL = {Bulletin of the American Mathematical Society},
    VOLUME = {83},
      YEAR = {1977},
    NUMBER = {5},
     PAGES = {831--879},
      ISSN = {0002-9904},
   MRCLASS = {55D35 (58A99)},
  MRNUMBER = {0454968 (56 \#13210)},
MRREVIEWER = {Jean-Michel Lemaire},
}

@book {Del2,
    AUTHOR = {Deligne, Pierre},
     TITLE = {\'{E}quations diff\'erentielles \`a points singuliers
              r\'eguliers},
    SERIES = {Lecture Notes in Mathematics, Vol. 163},
 PUBLISHER = {Springer-Verlag, Berlin-New York},
      YEAR = {1970},
     PAGES = {iii+133},
   MRCLASS = {14D05 (14C30)},
  MRNUMBER = {0417174},
MRREVIEWER = {Helmut Hamm},
}

@incollection {Deligne:P1,
   AUTHOR = {Deligne, P.},
    TITLE = {Le groupe fondamental de la droite projective moins trois
             points},
BOOKTITLE = {Galois groups over ${\bf Q}$ ({B}erkeley, {CA}, 1987)},
   SERIES = {Math. Sci. Res. Inst. Publ.},
   VOLUME = {16},
    PAGES = {79--297},
PUBLISHER = {Springer, New York},
     YEAR = {1989},
  MRCLASS = {14G25 (11G35 11M06 11R70 14F35 19E99 19F27)},
 MRNUMBER = {1012168 (90m:14016)},
MRREVIEWER = {James Milne},
      DOI = {10.1007/978-1-4613-9649-9_3},
      URL = {http://dx.doi.org/10.1007/978-1-4613-9649-9\_3},
}

@article {Deligne:Multizetas,
	AUTHOR = {Deligne, Pierre},
	TITLE = {Multiz\^etas, d'apr\`es {F}rancis {B}rown},
	NOTE = {S\'eminaire Bourbaki. Vol. 2011/2012. Expos\'es 1043--1058},
	JOURNAL = {Ast\'erisque},
	FJOURNAL = {Ast\'erisque},
	NUMBER = {352},
	YEAR = {2013},
	PAGES = {Exp. No. 1048, viii, 161--185},
	ISSN = {0303-1179},
	ISBN = {978-2-85629-371-3},
	MRCLASS = {11S40 (11G09 14C15 14F35)},
	MRNUMBER = {3087346},
	MRREVIEWER = {Damian R\~A\P ssler},
}

@article {DG,
    AUTHOR = {Deligne, Pierre and Goncharov, Alexander B.},
     TITLE = {Groupes fondamentaux motiviques de {T}ate mixte},
   JOURNAL = {Ann. Sci. \'Ecole Norm. Sup. (4)},
  FJOURNAL = {Annales Scientifiques de l'\'Ecole Normale Sup\'erieure.
              Quatri\`eme S\'erie},
    VOLUME = {38},
      YEAR = {2005},
    NUMBER = {1},
     PAGES = {1--56},
      ISSN = {0012-9593},
     CODEN = {ASENAH},
   MRCLASS = {11G55 (14F42 14G10 19F27)},
  MRNUMBER = {2136480 (2006b:11066)},
MRREVIEWER = {Tam{\'a}s Szamuely},
       DOI = {10.1016/j.ansens.2004.11.001},
       URL = {http://dx.doi.org/10.1016/j.ansens.2004.11.001},
}

@book {DM,
    AUTHOR = {Deligne, Pierre and Milne, James S. and Ogus, Arthur and Shih,
              Kuang-yen},
     TITLE = {Hodge cycles, motives, and {S}himura varieties},
    SERIES = {Lecture Notes in Mathematics},
    VOLUME = {900},
 PUBLISHER = {Springer-Verlag, Berlin-New York},
      YEAR = {1982},
     PAGES = {ii+414},
      ISBN = {3-540-11174-3},
   MRCLASS = {14Kxx (10D25 12A67 14A20 14F30 14K22)},
  MRNUMBER = {654325},
}

@preamble{
   "\def\cprime{$'$} "
}
@article {Drinfeld:Gal,
    AUTHOR = {Drinfel{\cprime}d, V. G.},
     TITLE = {On quasitriangular quasi-{H}opf algebras and on a group that
              is closely connected with {${\rm Gal}(\overline{\bf Q}/{\bf
              Q})$}},
   JOURNAL = {Algebra i Analiz},
  FJOURNAL = {Algebra i Analiz},
    VOLUME = {2},
      YEAR = {1990},
    NUMBER = {4},
     PAGES = {149--181},
      ISSN = {0234-0852},
   MRCLASS = {16W30 (17B37)},
  MRNUMBER = {1080203 (92f:16047)},
MRREVIEWER = {Ivan Penkov},
}

@preamble{
	"\def\cprime{$'$} "
}
@article {Drinfeld:QuasiHopf,
	AUTHOR = {Drinfel{\cprime}d, V. G.},
	TITLE = {Quasi-{H}opf algebras},
	JOURNAL = {Algebra i Analiz},
	FJOURNAL = {Algebra i Analiz},
	VOLUME = {1},
	YEAR = {1989},
	NUMBER = {6},
	PAGES = {114--148},
	ISSN = {0234-0852},
	MRCLASS = {17B37 (16W30 57M25 81T40)},
	MRNUMBER = {1047964},
	MRREVIEWER = {Ya. S. So{\u\i}bel{\cprime}man},
}

@article {Ecalle,
    AUTHOR = {Ecalle, Jean},
     TITLE = {A{RI}/{GARI}, la dimorphie et l'arithm\'etique des
              multiz\^etas: un premier bilan},
   JOURNAL = {J. Th\'eor. Nombres Bordeaux},
  FJOURNAL = {Journal de Th\'eorie des Nombres de Bordeaux},
    VOLUME = {15},
      YEAR = {2003},
    NUMBER = {2},
     PAGES = {411--478},
      ISSN = {1246-7405},
   MRCLASS = {11M41 (11G55 19F27 33B30)},
  MRNUMBER = {2140864},
MRREVIEWER = {Alexey A. Panchishkin},
       URL = {http://jtnb.cedram.org/item?id=JTNB_2003__15_2_411_0},
}

@article {Enriquez:EllAss,
    AUTHOR = {Enriquez, Benjamin},
     TITLE = {Elliptic associators},
   JOURNAL = {Selecta Math. (N.S.)},
  FJOURNAL = {Selecta Mathematica. New Series},
    VOLUME = {20},
      YEAR = {2014},
    NUMBER = {2},
     PAGES = {491--584},
      ISSN = {1022-1824},
   MRCLASS = {17B35 (11M32 14H10 16S30 20F36)},
  MRNUMBER = {3177926},
       DOI = {10.1007/s00029-013-0137-3},
       URL = {http://dx.doi.org/10.1007/s00029-013-0137-3},
}

@article {Enriquez:Emzv,
	AUTHOR = {Enriquez, Benjamin},
	TITLE = {Analogues elliptiques des nombres multiz\'etas},
	JOURNAL = {Bull. Soc. Math. France},
	FJOURNAL = {Bulletin de la Soci\'et\'e Math\'ematique de France},
	VOLUME = {144},
	YEAR = {2016},
	NUMBER = {3},
	PAGES = {395--427},
	ISSN = {0037-9484},
	MRCLASS = {11M32 (17B01 17B35 17B40 33E30)},
	MRNUMBER = {3558428},
	MRREVIEWER = {A. Perelli},
}

@article {Eichler,
	AUTHOR = {Eichler, M.},
	TITLE = {Eine {V}erallgemeinerung der {A}belschen {I}ntegrale},
	JOURNAL = {Math. Z.},
	FJOURNAL = {Mathematische Zeitschrift},
	VOLUME = {67},
	YEAR = {1957},
	PAGES = {267--298},
	ISSN = {0025-5874},
	MRCLASS = {33.0X},
	MRNUMBER = {0089928},
	MRREVIEWER = {H. Cohn},
}

@incollection {Fal,
    AUTHOR = {Faltings, Gerd},
     TITLE = {Mathematics around {K}im's new proof of {S}iegel's theorem},
 BOOKTITLE = {Diophantine geometry},
    SERIES = {CRM Series},
    VOLUME = {4},
     PAGES = {173--188},
 PUBLISHER = {Ed. Norm., Pisa},
      YEAR = {2007},
   MRCLASS = {14G99 (11G30 14F20)},
  MRNUMBER = {2349654 (2009i:14029)},
}

@article {FurStab,
    AUTHOR = {Furusho, Hidekazu},
     TITLE = {The multiple zeta value algebra and the stable derivation
              algebra},
   JOURNAL = {Publ. Res. Inst. Math. Sci.},
  FJOURNAL = {Kyoto University. Research Institute for Mathematical
              Sciences. Publications},
    VOLUME = {39},
      YEAR = {2003},
    NUMBER = {4},
     PAGES = {695--720},
      ISSN = {0034-5318},
     CODEN = {KRMPBV},
   MRCLASS = {11M41 (14G32)},
  MRNUMBER = {2025460},
       URL = {http://projecteuclid.org/euclid.prims/1145476044},
}

@incollection {FurMZVGT,
    AUTHOR = {Furusho, Hidekazu},
     TITLE = {Multiple zeta values and {G}rothendieck-{T}eichm\"uller
              groups},
 BOOKTITLE = {Primes and knots},
    SERIES = {Contemp. Math.},
    VOLUME = {416},
     PAGES = {49--82},
 PUBLISHER = {Amer. Math. Soc., Providence, RI},
      YEAR = {2006},
   MRCLASS = {14G32 (11M41)},
  MRNUMBER = {2276136},
MRREVIEWER = {Alexey A. Panchishkin},
       DOI = {10.1090/conm/416/07887},
       URL = {http://dx.doi.org/10.1090/conm/416/07887},
}

@article {Fur,
    AUTHOR = {Furusho, Hidekazu},
     TITLE = {Double shuffle relation for associators},
   JOURNAL = {Ann. of Math. (2)},
  FJOURNAL = {Annals of Mathematics. Second Series},
    VOLUME = {174},
      YEAR = {2011},
    NUMBER = {1},
     PAGES = {341--360},
      ISSN = {0003-486X},
     CODEN = {ANMAAH},
   MRCLASS = {14G32 (11G55 11M32 16W60)},
  MRNUMBER = {2811601 (2012i:14031)},
MRREVIEWER = {Pierre A. Lochak},
       DOI = {10.4007/annals.2011.174.1.9},
       URL = {http://dx.doi.org/10.4007/annals.2011.174.1.9},
}

@incollection {GKZ,
    AUTHOR = {Gangl, Herbert and Kaneko, Masanobu and Zagier, Don},
     TITLE = {Double zeta values and modular forms},
 BOOKTITLE = {Automorphic forms and zeta functions},
     PAGES = {71--106},
 PUBLISHER = {World Sci. Publ., Hackensack, NJ},
      YEAR = {2006},
   MRCLASS = {11M41 (11F11)},
  MRNUMBER = {2208210 (2006m:11138)},
MRREVIEWER = {Hirofumi Tsumura},
       DOI = {10.1142/9789812774415_0004},
       URL = {http://dx.doi.org/10.1142/9789812774415_0004},
}

@article {GonMod,
    AUTHOR = {Goncharov, A. B.},
     TITLE = {Multiple polylogarithms, cyclotomy and modular complexes},
   JOURNAL = {Math. Res. Lett.},
  FJOURNAL = {Mathematical Research Letters},
    VOLUME = {5},
      YEAR = {1998},
    NUMBER = {4},
     PAGES = {497--516},
      ISSN = {1073-2780},
   MRCLASS = {11G55 (11F67 11R42 19E20 19F15 19F27)},
  MRNUMBER = {1653320 (2000c:11108)},
MRREVIEWER = {Alexey A. Panchishkin},
       DOI = {10.4310/MRL.1998.v5.n4.a7},
       URL = {http://dx.doi.org/10.4310/MRL.1998.v5.n4.a7},
}

@incollection {GonMZV,
    AUTHOR = {Goncharov, Alexander B.},
     TITLE = {Multiple {$\zeta$}-values, {G}alois groups, and geometry of
              modular varieties},
 BOOKTITLE = {European {C}ongress of {M}athematics, {V}ol. {I} ({B}arcelona,
              2000)},
    SERIES = {Progr. Math.},
    VOLUME = {201},
     PAGES = {361--392},
 PUBLISHER = {Birkh\"auser, Basel},
      YEAR = {2001},
   MRCLASS = {11G55 (11G40 11M41 14G35 33B30 81Q30)},
  MRNUMBER = {1905330},
MRREVIEWER = {Jan Nekov{\'a}{\v{r}}},
}

@unpublished{Goncharov:MTM,
   author = {{Goncharov}, A.~B.},
    title = "{Multiple polylogarithms and mixed Tate motives}",
    note = {arXiv:math/0103059},
    year = 2001,
}

@article {GM,
    AUTHOR = {Goncharov, A. B. and Manin, Yu. I.},
     TITLE = {Multiple {$\zeta$}-motives and moduli spaces {$\overline{\mathscr
              M}_{0,n}$}},
   JOURNAL = {Compos. Math.},
  FJOURNAL = {Compositio Mathematica},
    VOLUME = {140},
      YEAR = {2004},
    NUMBER = {1},
     PAGES = {1--14},
      ISSN = {0010-437X},
   MRCLASS = {11G55 (11M41 14H10)},
  MRNUMBER = {2004120 (2005c:11090)},
MRREVIEWER = {Gilberto Bini},
       DOI = {10.1112/S0010437X03000125},
       URL = {http://dx.doi.org/10.1112/S0010437X03000125},
}

@incollection {Hai,
    AUTHOR = {Hain, Richard M.},
     TITLE = {The geometry of the mixed {H}odge structure on the fundamental
              group},
 BOOKTITLE = {Algebraic geometry, {B}owdoin, 1985 ({B}runswick, {M}aine,
              1985)},
    SERIES = {Proc. Sympos. Pure Math.},
    VOLUME = {46},
     PAGES = {247--282},
 PUBLISHER = {Amer. Math. Soc., Providence, RI},
      YEAR = {1987},
   MRCLASS = {14F40 (14C30 32C40 32G20 55P62 58A14)},
  MRNUMBER = {927984 (89g:14010)},
MRREVIEWER = {Toshitake Kohno},
}

@incollection {HaiPolylog,
    AUTHOR = {Hain, Richard M.},
     TITLE = {Classical polylogarithms},
 BOOKTITLE = {Motives ({S}eattle, {WA}, 1991)},
    SERIES = {Proc. Sympos. Pure Math.},
    VOLUME = {55},
     PAGES = {3--42},
 PUBLISHER = {Amer. Math. Soc., Providence, RI},
      YEAR = {1994},
   MRCLASS = {19F99 (11G99 11R70 19D55)},
  MRNUMBER = {1265550 (94k:19002)},
MRREVIEWER = {Philippe Blanc},
}

@unpublished{Hain:KZB,
	author = {{Hain}, R.},
	title = "{Notes on the Universal Elliptic KZB Equation}",
	note = {arXiv:1309.0580},
	year = 2013,
}

@incollection {Hain:HodgeDeRham,
	AUTHOR = {Hain, Richard},
	TITLE = {The {H}odge--de {R}ham theory of modular groups},
	BOOKTITLE = {Recent advances in {H}odge theory},
	SERIES = {London Math. Soc. Lecture Note Ser.},
	VOLUME = {427},
	PAGES = {422--514},
	PUBLISHER = {Cambridge Univ. Press, Cambridge},
	YEAR = {2016},
	MRCLASS = {14D07 (14Gxx 32S35 58A14)},
	MRNUMBER = {3409885},
}

@unpublished{HM,
   author = {{Hain}, R. and {Matsumoto}, M.},
    title = "{Universal Mixed Elliptic Motives}",
  note = {arXiv:1512.03975},
     year = 2015
}

@article {Ihara:Annals,
    AUTHOR = {Ihara, Yasutaka},
     TITLE = {Profinite braid groups, {G}alois representations and complex
              multiplications},
   JOURNAL = {Ann. of Math. (2)},
  FJOURNAL = {Annals of Mathematics. Second Series},
    VOLUME = {123},
      YEAR = {1986},
    NUMBER = {1},
     PAGES = {43--106},
      ISSN = {0003-486X},
     CODEN = {ANMAAH},
   MRCLASS = {11G25 (14K22 20E18)},
  MRNUMBER = {825839},
MRREVIEWER = {David Goss},
       DOI = {10.2307/1971352},
       URL = {http://dx.doi.org/10.2307/1971352},
}

@incollection {Iha,
    AUTHOR = {Ihara, Yasutaka},
     TITLE = {The {G}alois representation arising from {${\bf P}^1-\{0,1,\infty\}$} and {T}ate twists of even degree},
 BOOKTITLE = {Galois groups over {${\bf Q}$} ({B}erkeley, {CA}, 1987)},
    SERIES = {Math. Sci. Res. Inst. Publ.},
    VOLUME = {16},
     PAGES = {299--313},
 PUBLISHER = {Springer, New York},
      YEAR = {1989},
   MRCLASS = {11F80 (11G20 11R23 11R58 14E22 14G25)},
  MRNUMBER = {1012169},
MRREVIEWER = {Sheldon Kamienny},
       DOI = {10.1007/978-1-4613-9649-9_4},
       URL = {http://dx.doi.org/10.1007/978-1-4613-9649-9_4},
}

@inproceedings {IhICM,
    AUTHOR = {Ihara, Yasutaka},
     TITLE = {Braids, {G}alois groups, and some arithmetic functions},
 BOOKTITLE = {Proceedings of the {I}nternational {C}ongress of
              {M}athematicians, {V}ol.\ {I}, {II} ({K}yoto, 1990)},
     PAGES = {99--120},
 PUBLISHER = {Math. Soc. Japan, Tokyo},
      YEAR = {1991},
   MRCLASS = {11G09 (11R32 14E20 16W30 20F34)},
  MRNUMBER = {1159208},
MRREVIEWER = {J. Browkin},
}

@incollection {IKY,
    AUTHOR = {Ihara, Yasutaka and Kaneko, Masanobu and Yukinari, Atsushi},
     TITLE = {On some properties of the universal power series for {J}acobi
              sums},
 BOOKTITLE = {Galois representations and arithmetic algebraic geometry
              ({K}yoto, 1985/{T}okyo, 1986)},
    SERIES = {Adv. Stud. Pure Math.},
    VOLUME = {12},
     PAGES = {65--86},
 PUBLISHER = {North-Holland, Amsterdam},
      YEAR = {1987},
   MRCLASS = {11R23 (11S80)},
  MRNUMBER = {948237},
MRREVIEWER = {J. Browkin},
}

@article {IKZ,
    AUTHOR = {Ihara, Kentaro and Kaneko, Masanobu and Zagier, Don},
     TITLE = {Derivation and double shuffle relations for multiple zeta
              values},
   JOURNAL = {Compos. Math.},
  FJOURNAL = {Compositio Mathematica},
    VOLUME = {142},
      YEAR = {2006},
    NUMBER = {2},
     PAGES = {307--338},
      ISSN = {0010-437X},
   MRCLASS = {11M41},
  MRNUMBER = {2218898},
MRREVIEWER = {David Bradley},
       DOI = {10.1112/S0010437X0500182X},
       URL = {http://dx.doi.org/10.1112/S0010437X0500182X},
}

@article {IO,
    AUTHOR = {Ihara, Kentaro and Ochiai, Hiroyuki},
     TITLE = {Symmetry on linear relations for multiple zeta values},
   JOURNAL = {Nagoya Math. J.},
  FJOURNAL = {Nagoya Mathematical Journal},
    VOLUME = {189},
      YEAR = {2008},
     PAGES = {49--62},
      ISSN = {0027-7630},
     CODEN = {NGMJA2},
   MRCLASS = {11M41},
  MRNUMBER = {2396583},
MRREVIEWER = {David Bradley},
       URL = {http://projecteuclid.org/euclid.nmj/1205156910},
}

@book {Kas,
    AUTHOR = {Kassel, Christian},
     TITLE = {Quantum groups},
    SERIES = {Graduate Texts in Mathematics},
    VOLUME = {155},
 PUBLISHER = {Springer-Verlag, New York},
      YEAR = {1995},
     PAGES = {xii+531},
      ISBN = {0-387-94370-6},
   MRCLASS = {17B37 (16W30 18D10 20F36 57M25 81R50)},
  MRNUMBER = {1321145 (96e:17041)},
MRREVIEWER = {Yu. N. Bespalov},
       DOI = {10.1007/978-1-4612-0783-2},
       URL = {http://dx.doi.org/10.1007/978-1-4612-0783-2},
}

@article {Kim:Annals,
    AUTHOR = {Kim, Minhyong},
     TITLE = {{$p$}-adic {$L$}-functions and {S}elmer varieties associated
              to elliptic curves with complex multiplication},
   JOURNAL = {Ann. of Math. (2)},
  FJOURNAL = {Annals of Mathematics. Second Series},
    VOLUME = {172},
      YEAR = {2010},
    NUMBER = {1},
     PAGES = {751--759},
      ISSN = {0003-486X},
     CODEN = {ANMAAH},
   MRCLASS = {11G15 (11G05 11G40)},
  MRNUMBER = {2680431 (2011i:11089)},
MRREVIEWER = {Francesc C. Castell{\`a}},
       DOI = {10.4007/annals.2010.172.751},
       URL = {http://dx.doi.org/10.4007/annals.2010.172.751},
}

@article {KZ,
    AUTHOR = {Knizhnik, V. G. and Zamolodchikov, A. B.},
     TITLE = {Current algebra and {W}ess-{Z}umino model in two dimensions},
   JOURNAL = {Nuclear Phys. B},
  FJOURNAL = {Nuclear Physics. B},
    VOLUME = {247},
      YEAR = {1984},
    NUMBER = {1},
     PAGES = {83--103},
      ISSN = {0550-3213},
     CODEN = {NUPBBO},
   MRCLASS = {81E13 (81D15)},
  MRNUMBER = {853258},
       DOI = {10.1016/0550-3213(84)90374-2},
       URL = {http://dx.doi.org/10.1016/0550-3213(84)90374-2},
}

@book {Lang,
    AUTHOR = {Lang, Serge},
     TITLE = {Introduction to modular forms},
      NOTE = {Grundlehren der mathematischen Wissenschaften, No. 222},
 PUBLISHER = {Springer-Verlag, Berlin-New York},
      YEAR = {1976},
     PAGES = {ix+261},
   MRCLASS = {10DXX},
  MRNUMBER = {0429740},
MRREVIEWER = {Neal Koblitz},
}

@article {Landen,
    AUTHOR = {Landen, John},
     TITLE = {Mathematical memoirs respecting a variety of subjects},
   JOURNAL = {Nourse},
      YEAR = {1780},
} 

@article {LM,
    AUTHOR = {Le, Thang Tu Quoc and Murakami, Jun},
     TITLE = {Kontsevich's integral for the {K}auffman polynomial},
   JOURNAL = {Nagoya Math. J.},
  FJOURNAL = {Nagoya Mathematical Journal},
    VOLUME = {142},
      YEAR = {1996},
     PAGES = {39--65},
      ISSN = {0027-7630},
     CODEN = {NGMJA2},
   MRCLASS = {57M25 (11M99)},
  MRNUMBER = {1399467 (97d:57009)},
MRREVIEWER = {Sergei K. Lando},
       URL = {http://projecteuclid.org/euclid.nmj/1118772043},
}

@article {LM:Compositio,
	AUTHOR = {Le, Tu Quoc Thang and Murakami, Jun},
	TITLE = {The universal {V}assiliev-{K}ontsevich invariant for framed
		oriented links},
	JOURNAL = {Compositio Math.},
	FJOURNAL = {Compositio Mathematica},
	VOLUME = {102},
	YEAR = {1996},
	NUMBER = {1},
	PAGES = {41--64},
	ISSN = {0010-437X},
	CODEN = {CMPMAF},
	MRCLASS = {57M25 (11M41)},
	MRNUMBER = {1394520},
	URL = {http://www.numdam.org/item?id=CM_1996__102_1_41_0},
}
 
@article {Levin:Compositio,
    AUTHOR = {Levin, Andrey},
     TITLE = {Elliptic polylogarithms: an analytic theory},
   JOURNAL = {Compositio Math.},
  FJOURNAL = {Compositio Mathematica},
    VOLUME = {106},
      YEAR = {1997},
    NUMBER = {3},
     PAGES = {267--282},
      ISSN = {0010-437X},
     CODEN = {CMPMAF},
   MRCLASS = {11F37 (11F27 11G40 11R70 19F27)},
  MRNUMBER = {1457106 (98d:11048)},
MRREVIEWER = {Alexey A. Panchishkin},
       DOI = {10.1023/A:1000193320513},
       URL = {http://dx.doi.org/10.1023/A:1000193320513},
}

@unpublished{LR,
     author         = {A. Levin and G. Racinet},
     title          = {Towards multiple elliptic polylogarithms},
	 note   		= {arXiv:math/0703237},
     year           = {2007},
}

@unpublished{LMS,
	author = {Lochak, Pierre and Matthes, Nils and Schneps, Leila},
	title = "{Elliptic multiple zeta values and the elliptic double shuffle relations}",
	year = {2017},
	note = {arXiv:1703.09410},
}

@incollection {Manin:Iterated,
    AUTHOR = {Manin, Yuri I.},
     TITLE = {Iterated integrals of modular forms and noncommutative modular
              symbols},
 BOOKTITLE = {Algebraic geometry and number theory},
    SERIES = {Progr. Math.},
    VOLUME = {253},
     PAGES = {565--597},
 PUBLISHER = {Birkh\"auser Boston, Boston, MA},
      YEAR = {2006},
   MRCLASS = {11F67 (11G55 11M41)},
  MRNUMBER = {2263200 (2008a:11062)},
MRREVIEWER = {Caterina Consani},
       DOI = {10.1007/978-0-8176-4532-8_10},
       URL = {http://dx.doi.org/10.1007/978-0-8176-4532-8\_10},
}

@article {Matthes:Edzv,
	AUTHOR = {Matthes, Nils},
	TITLE = {Elliptic double zeta values},
	JOURNAL = {J. Number Theory},
	FJOURNAL = {Journal of Number Theory},
	VOLUME = {171},
	YEAR = {2017},
	PAGES = {227--251},
	ISSN = {0022-314X},
	CODEN = {JNUTA9},
	MRCLASS = {11M32 (11F50)},
	MRNUMBER = {3556684},
	DOI = {10.1016/j.jnt.2016.07.010},
	URL = {http://dx.doi.org/10.1016/j.jnt.2016.07.010},
}

@PhDThesis{Matthes:Thesis,
	author     =     {Matthes, Nils},
	title     =     {{Elliptic multiple zeta values}},
	school     =     {Universit\"at Hamburg},
	year     =     {2016},
}

@unpublished{Matthes:Metab,
	author = {{Matthes}, N.},
	title = "{The meta-abelian elliptic KZB associator and periods of Eisenstein series}",
	note = {arXiv:1608.00740},
	year = 2016,
}

@article {Nakamura:Galoisrep,
    AUTHOR = {Nakamura, Hiroaki},
     TITLE = {On exterior {G}alois representations associated with open
              elliptic curves},
   JOURNAL = {J. Math. Sci. Univ. Tokyo},
  FJOURNAL = {The University of Tokyo. Journal of Mathematical Sciences},
    VOLUME = {2},
      YEAR = {1995},
    NUMBER = {1},
     PAGES = {197--231},
      ISSN = {1340-5705},
   MRCLASS = {11G05 (11F80 11G16 14H30)},
  MRNUMBER = {1348028},
MRREVIEWER = {Yasutaka Ihara},
}

@incollection {Nakamura:Eisen99,
    AUTHOR = {Nakamura, Hiroaki},
     TITLE = {Tangential base points and {E}isenstein power series},
 BOOKTITLE = {Aspects of {G}alois theory ({G}ainesville, {FL}, 1996)},
    SERIES = {London Math. Soc. Lecture Note Ser.},
    VOLUME = {256},
     PAGES = {202--217},
 PUBLISHER = {Cambridge Univ. Press, Cambridge},
      YEAR = {1999},
   MRCLASS = {14G32 (11G07 11G55)},
  MRNUMBER = {1708607},
MRREVIEWER = {Helmut V{\"o}lklein},
}

@article {Nakamura:Arithmetic,
	AUTHOR = {Nakamura, Hiroaki},
	TITLE = {On arithmetic monodromy representations of {E}isenstein type
		in fundamental groups of once punctured elliptic curves},
	JOURNAL = {Publ. Res. Inst. Math. Sci.},
	FJOURNAL = {Publications of the Research Institute for Mathematical
		Sciences},
	VOLUME = {49},
	YEAR = {2013},
	NUMBER = {3},
	PAGES = {413--496},
	ISSN = {0034-5318},
	MRCLASS = {14G32 (11F20 11G16 14G25)},
	MRNUMBER = {3097013},
	MRREVIEWER = {Kirsten Wickelgren},
	DOI = {10.4171/PRIMS/110},
	URL = {http://dx.doi.org/10.4171/PRIMS/110},
}

@unpublished{Nakamura:EisenRevisited,
	Author = {Nakamura,Hiroaki},
	Title = {On profinite {E}isenstein periods in the monodromy of universal elliptic curves},
	YEAR = {2016},
	note = {http://www.math.sci.osaka-u.ac.jp/$\sim$nakamura/zoo/fox/EisenRevisited.pdf},
}	

@MastersThesis{Pollack:Thesis,
    author     =     {Pollack, Aaron},
    title     =     {{Relations between derivations arising from modular forms}},
    school     =     {Duke University},
    year     =     {2009},
    }

@article {Rac,
    AUTHOR = {Racinet, Georges},
     TITLE = {Doubles m\'elanges des polylogarithmes multiples aux racines
              de l'unit\'e},
   JOURNAL = {Publ. Math. Inst. Hautes \'Etudes Sci.},
  FJOURNAL = {Publications Math\'ematiques. Institut de Hautes \'Etudes
              Scientifiques},
    NUMBER = {95},
      YEAR = {2002},
     PAGES = {185--231},
      ISSN = {0073-8301},
   MRCLASS = {11G55 (11M41)},
  MRNUMBER = {1953193 (2004c:11117)},
MRREVIEWER = {Jan Nekov{\'a}{\v{r}}},
       DOI = {10.1007/s102400200004},
       URL = {http://dx.doi.org/10.1007/s102400200004},
}

@article {Ree,
    AUTHOR = {Ree, Rimhak},
     TITLE = {Lie elements and an algebra associated with shuffles},
   JOURNAL = {Ann. of Math. (2)},
  FJOURNAL = {Annals of Mathematics. Second Series},
    VOLUME = {68},
      YEAR = {1958},
     PAGES = {210--220},
      ISSN = {0003-486X},
   MRCLASS = {17.00 (20.00)},
  MRNUMBER = {0100011 (20 \#6447)},
MRREVIEWER = {P. M. Cohn},
}

@book {Reutenauer,
    AUTHOR = {Reutenauer, Christophe},
     TITLE = {Free {L}ie algebras},
    SERIES = {London Mathematical Society Monographs. New Series},
    VOLUME = {7},
      NOTE = {Oxford Science Publications},
 PUBLISHER = {The Clarendon Press, Oxford University Press, New York},
      YEAR = {1993},
     PAGES = {xviii+269},
      ISBN = {0-19-853679-8},
   MRCLASS = {17-02 (05-02 17B05)},
  MRNUMBER = {1231799 (94j:17002)},
MRREVIEWER = {Hartmut Laue},
}

@article {SS,
    AUTHOR = {Schlotterer, O. and Stieberger, S.},
     TITLE = {Motivic multiple zeta values and superstring amplitudes},
   JOURNAL = {J. Phys. A},
  FJOURNAL = {Journal of Physics. A. Mathematical and Theoretical},
    VOLUME = {46},
      YEAR = {2013},
    NUMBER = {47},
     PAGES = {475401, 37},
      ISSN = {1751-8113},
   MRCLASS = {81T30 (14E18)},
  MRNUMBER = {3126883},
MRREVIEWER = {Jihye Sofia Seo},
       DOI = {10.1088/1751-8113/46/47/475401},
       URL = {http://dx.doi.org/10.1088/1751-8113/46/47/475401},
}

@ARTICLE{Sch,
   author = {{Schneps}, L.},
    title = "{Elliptic multiple zeta values, Grothendieck-Teichm{\"u}ller and mould theory}",
  journal = {ArXiv e-prints},
   volume = {math.NT/1506.09050},
   year = 2015
}

@book {SerArith,
    AUTHOR = {Serre, Jean-Pierre},
     TITLE = {Cours d'arithm\'etique},
    SERIES = {Collection SUP: ``Le Math\'ematicien''},
    VOLUME = {2},
 PUBLISHER = {Presses Universitaires de France, Paris},
      YEAR = {1970},
     PAGES = {188},
   MRCLASS = {10.01},
  MRNUMBER = {0255476},
MRREVIEWER = {Burton W. Jones},
}

@book {Ser,
    AUTHOR = {Serre, Jean-Pierre},
     TITLE = {Repr\'esentations lin\'eaires des groupes finis},
   EDITION = {revised},
 PUBLISHER = {Hermann, Paris},
      YEAR = {1978},
     PAGES = {182},
      ISBN = {2-7056-5630-8},
   MRCLASS = {20-01 (20C99)},
  MRNUMBER = {543841 (80f:20001)},
}

@book {Serre:Lie,
    AUTHOR = {Serre, Jean-Pierre},
     TITLE = {Lie algebras and {L}ie groups},
    SERIES = {Lecture Notes in Mathematics},
    VOLUME = {1500},
      NOTE = {1964 lectures given at Harvard University,
              Corrected fifth printing of the second (1992) edition},
 PUBLISHER = {Springer-Verlag, Berlin},
      YEAR = {2006},
     PAGES = {viii+168},
      ISBN = {978-3-540-55008-2; 3-540-55008-9},
   MRCLASS = {17-01 (22-01)},
  MRNUMBER = {2179691},
}

@incollection {Sou,
    AUTHOR = {Soul{\'e}, Christophe},
     TITLE = {On higher {$p$}-adic regulators},
 BOOKTITLE = {Algebraic {$K$}-theory, {E}vanston 1980 ({P}roc. {C}onf.,
              {N}orthwestern {U}niv., {E}vanston, {I}ll., 1980)},
    SERIES = {Lecture Notes in Math.},
    VOLUME = {854},
     PAGES = {372--401},
 PUBLISHER = {Springer, Berlin-New York},
      YEAR = {1981},
   MRCLASS = {12A62 (12B22 13D03 18F25)},
  MRNUMBER = {618313},
MRREVIEWER = {J. Browkin},
}

@article {Terasoma:MTM,
    AUTHOR = {Terasoma, Tomohide},
     TITLE = {Mixed {T}ate motives and multiple zeta values},
   JOURNAL = {Invent. Math.},
  FJOURNAL = {Inventiones Mathematicae},
    VOLUME = {149},
      YEAR = {2002},
    NUMBER = {2},
     PAGES = {339--369},
      ISSN = {0020-9910},
     CODEN = {INVMBH},
   MRCLASS = {11G55 (11M41 19F27)},
  MRNUMBER = {1918675},
MRREVIEWER = {Jan Nekov{\'a}{\v{r}}},
       DOI = {10.1007/s002220200218},
       URL = {http://dx.doi.org/10.1007/s002220200218},
}

@incollection {Tera,
    AUTHOR = {Terasoma, Tomohide},
     TITLE = {Geometry of multiple zeta values},
 BOOKTITLE = {International {C}ongress of {M}athematicians. {V}ol. {II}},
     PAGES = {627--635},
 PUBLISHER = {Eur. Math. Soc., Z\"urich},
      YEAR = {2006},
   MRCLASS = {14C30 (11M41 14F42)},
  MRNUMBER = {2275614 (2008e:14009)},
MRREVIEWER = {Jan Nekov{\'a}{\v{r}}},
}

@article {Tsunogai:Derivations,
    AUTHOR = {Tsunogai, Hiroshi},
     TITLE = {On some derivations of {L}ie algebras related to {G}alois
              representations},
   JOURNAL = {Publ. Res. Inst. Math. Sci.},
  FJOURNAL = {Kyoto University. Research Institute for Mathematical
              Sciences. Publications},
    VOLUME = {31},
      YEAR = {1995},
    NUMBER = {1},
     PAGES = {113--134},
      ISSN = {0034-5318},
     CODEN = {KRMPBV},
   MRCLASS = {11G05 (11R32 14H30 17B40)},
  MRNUMBER = {1317526},
MRREVIEWER = {Douglas L. Ulmer},
       DOI = {10.2977/prims/1195164794},
       URL = {http://dx.doi.org/10.2977/prims/1195164794},
}

@article {Tsumu,
    AUTHOR = {Tsumura, Hirofumi},
     TITLE = {Combinatorial relations for {E}uler-{Z}agier sums},
   JOURNAL = {Acta Arith.},
  FJOURNAL = {Acta Arithmetica},
    VOLUME = {111},
      YEAR = {2004},
    NUMBER = {1},
     PAGES = {27--42},
      ISSN = {0065-1036},
     CODEN = {AARIA9},
   MRCLASS = {11M41 (33E20)},
  MRNUMBER = {2038060},
MRREVIEWER = {David Bradley},
       DOI = {10.4064/aa111-1-3},
       URL = {http://dx.doi.org/10.4064/aa111-1-3},
}

@book {Was,
    AUTHOR = {Wasow, Wolfgang},
     TITLE = {Asymptotic expansions for ordinary differential equations},
      NOTE = {Reprint of the 1976 edition},
 PUBLISHER = {Dover Publications, Inc., New York},
      YEAR = {1987},
     PAGES = {x+374},
      ISBN = {0-486-65456-7},
   MRCLASS = {34-02 (34E05)},
  MRNUMBER = {919406},
}

@book {Wat,
    AUTHOR = {Waterhouse, William C.},
     TITLE = {Introduction to affine group schemes},
    SERIES = {Graduate Texts in Mathematics},
    VOLUME = {66},
 PUBLISHER = {Springer-Verlag, New York-Berlin},
      YEAR = {1979},
     PAGES = {xi+164},
      ISBN = {0-387-90421-2},
   MRCLASS = {14-01 (14Lxx 20G99)},
  MRNUMBER = {547117 (82e:14003)},
MRREVIEWER = {M. Kh. Gizatullin},
}

@book {W,
    AUTHOR = {Weil, Andr{\'e}},
     TITLE = {Elliptic functions according to {E}isenstein and {K}ronecker},
      NOTE = {Ergebnisse der Mathematik und ihrer Grenzgebiete, Band 88},
 PUBLISHER = {Springer-Verlag, Berlin-New York},
      YEAR = {1976},
     PAGES = {ii+93},
      ISBN = {3-540-07422-8},
   MRCLASS = {10DXX (01A55 10-03)},
  MRNUMBER = {0562289 (58 \#27769a)},
MRREVIEWER = {S. Chowla},
}

@article {ZagEll,
    AUTHOR = {Zagier, Don},
     TITLE = {The {B}loch-{W}igner-{R}amakrishnan polylogarithm function},
   JOURNAL = {Math. Ann.},
  FJOURNAL = {Mathematische Annalen},
    VOLUME = {286},
      YEAR = {1990},
    NUMBER = {1-3},
     PAGES = {613--624},
      ISSN = {0025-5831},
     CODEN = {MAANA},
   MRCLASS = {11R42 (11R70 19F27)},
  MRNUMBER = {1032949 (90k:11153)},
MRREVIEWER = {V. Kumar Murty},
       DOI = {10.1007/BF01453591},
       URL = {http://dx.doi.org/10.1007/BF01453591},
}

@article {Zagier:Periods,
    AUTHOR = {Zagier, Don},
     TITLE = {Periods of modular forms and {J}acobi theta functions},
   JOURNAL = {Invent. Math.},
  FJOURNAL = {Inventiones Mathematicae},
    VOLUME = {104},
      YEAR = {1991},
    NUMBER = {3},
     PAGES = {449--465},
      ISSN = {0020-9910},
     CODEN = {INVMBH},
   MRCLASS = {11F67 (11F27 11F55)},
  MRNUMBER = {1106744 (92e:11052)},
MRREVIEWER = {Rolf Berndt},
       DOI = {10.1007/BF01245085},
       URL = {http://dx.doi.org/10.1007/BF01245085},
}

@incollection {Zag,
    AUTHOR = {Zagier, Don},
     TITLE = {Values of zeta functions and their applications},
 BOOKTITLE = {First {E}uropean {C}ongress of {M}athematics, {V}ol.\ {II}
              ({P}aris, 1992)},
    SERIES = {Progr. Math.},
    VOLUME = {120},
     PAGES = {497--512},
 PUBLISHER = {Birkh\"auser, Basel},
      YEAR = {1994},
   MRCLASS = {11M41 (11F67 11G40 19F27)},
  MRNUMBER = {1341859 (96k:11110)},
MRREVIEWER = {Fernando Rodr{\'{\i}}guez Villegas},
}

@article {Zagier:Traces,
    AUTHOR = {Zagier, Don},
     TITLE = {Periods of modular forms, traces of {H}ecke operators, and
              multiple zeta values},
      NOTE = {Research into automorphic forms and $L$ functions (Japanese)
              (Kyoto, 1992)},
   JOURNAL = {S\=urikaisekikenky\=usho K\=oky\=uroku},
  FJOURNAL = {S\=urikaisekikenky\=usho K\=oky\=uroku},
    NUMBER = {843},
      YEAR = {1993},
     PAGES = {162--170},
   MRCLASS = {11F67 (11F25 11F72)},
  MRNUMBER = {1296720},
MRREVIEWER = {Pavel Guerzhoy},
}

@incollection {123,
    AUTHOR = {Zagier, Don},
     TITLE = {Elliptic modular forms and their applications},
 BOOKTITLE = {The 1-2-3 of modular forms},
    SERIES = {Universitext},
     PAGES = {1--103},
 PUBLISHER = {Springer, Berlin},
      YEAR = {2008},
   MRCLASS = {11F11 (11-02 11E45 11F20 11F25 11F27 11F67)},
  MRNUMBER = {2409678},
MRREVIEWER = {Rainer Schulze-Pillot},
       DOI = {10.1007/978-3-540-74119-0_1},
       URL = {http://dx.doi.org/10.1007/978-3-540-74119-0_1},
}
		
@article {Zud,
    AUTHOR = {Zudilin, V. V.},
     TITLE = {One of the numbers {$\zeta(5)$}, {$\zeta(7)$}, {$\zeta(9)$},
              {$\zeta(11)$} is irrational},
   JOURNAL = {Uspekhi Mat. Nauk},
  FJOURNAL = {Rossi\u\i skaya Akademiya Nauk. Moskovskoe Matematicheskoe
              Obshchestvo. Uspekhi Matematicheskikh Nauk},
    VOLUME = {56},
      YEAR = {2001},
    NUMBER = {4(340)},
     PAGES = {149--150},
      ISSN = {0042-1316},
   MRCLASS = {11J72 (11J91 11M06)},
  MRNUMBER = {1861452},
MRREVIEWER = {John H. Loxton},
       DOI = {10.1070/RM2001v056n04ABEH000427},
       URL = {http://dx.doi.org/10.1070/RM2001v056n04ABEH000427},
}

\end{bibtex}
\bibliographystyle{abbrv}
{\bibliography{\jobname}}

\def\cprime{$'$} \def\cprime{$'$}
\begin{thebibliography}{10}

\bibitem{BeiLev}
A.~Be{\u\i}linson and A.~Levin.
\newblock The elliptic polylogarithm.
\newblock In {\em Motives ({S}eattle, {WA}, 1991)}, volume~55 of {\em Proc.
  Sympos. Pure Math.}, pages 123--190. Amer. Math. Soc., Providence, RI, 1994.

\bibitem{BMS}
J.~Broedel, N.~Matthes, and O.~Schlotterer.
\newblock Relations between elliptic multiple zeta values and a special
  derivation algebra.
\newblock {\em J. Phys. A}, 49(15):155--203, 2016.

\bibitem{Brown:MTM}
F.~Brown.
\newblock Mixed {T}ate motives over {$\mathbb Z$}.
\newblock {\em Ann. of Math. (2)}, 175(2):949--976, 2012.

\bibitem{Brown:Decomposition}
F.~Brown.
\newblock On the decomposition of motivic multiple zeta values.
\newblock In {\em Galois-{T}eichm\"uller theory and arithmetic geometry},
  volume~63 of {\em Adv. Stud. Pure Math.}, pages 31--58. Math. Soc. Japan,
  Tokyo, 2012.

\bibitem{Brown:Depth}
F.~Brown.
\newblock Depth-graded motivic multiple zeta values.
\newblock arXiv:1301.3053, 2013.

\bibitem{Brown:MMV}
F.~Brown.
\newblock Multiple modular values and the relative completion of the
  fundamental group of $\mathcal{M}_{1,1}$.
\newblock arXiv:1407.5167v3, 2016.

\bibitem{Brown:Depth3}
F.~Brown.
\newblock Zeta elements in depth 3 and the fundamental {L}ie algebra of the
  infinitesimal {T}ate curve.
\newblock {\em Forum Math. Sigma}, 5, e1:56pp, 2017.

\bibitem{BL:MEP}
F.~Brown and A.~Levin.
\newblock Multiple elliptic polylogarithms.
\newblock arXiv:1110.6917, 2011.

\bibitem{CEE:KZB}
D.~Calaque, B.~Enriquez, and P.~Etingof.
\newblock Universal {KZB} equations: the elliptic case.
\newblock In {\em Algebra, arithmetic, and geometry: in honor of {Y}u. {I}.
  {M}anin. {V}ol. {I}}, volume 269 of {\em Progr. Math.}, pages 165--266.
  Birkh\"auser Boston, Inc., Boston, MA, 2009.

\bibitem{Chen:PathIntegrals}
K.~T. Chen.
\newblock Iterated path integrals.
\newblock {\em Bull. Amer. Math. Soc.}, 83(5):831--879, 1977.

\bibitem{Deligne:P1}
P.~Deligne.
\newblock Le groupe fondamental de la droite projective moins trois points.
\newblock In {\em Galois groups over ${\bf Q}$ ({B}erkeley, {CA}, 1987)},
  volume~16 of {\em Math. Sci. Res. Inst. Publ.}, pages 79--297. Springer, New
  York, 1989.

\bibitem{DG}
P.~Deligne and A.~B. Goncharov.
\newblock Groupes fondamentaux motiviques de {T}ate mixte.
\newblock {\em Ann. Sci. \'Ecole Norm. Sup. (4)}, 38(1):1--56, 2005.

\bibitem{Drinfeld:Gal}
V.~G. Drinfel{\cprime}d.
\newblock On quasitriangular quasi-{H}opf algebras and on a group that is
  closely connected with {${\rm Gal}(\overline{\bf Q}/{\bf Q})$}.
\newblock {\em Algebra i Analiz}, 2(4):149--181, 1990.

\bibitem{Enriquez:EllAss}
B.~Enriquez.
\newblock Elliptic associators.
\newblock {\em Selecta Math. (N.S.)}, 20(2):491--584, 2014.

\bibitem{Enriquez:Emzv}
B.~Enriquez.
\newblock Analogues elliptiques des nombres multiz\'etas.
\newblock {\em Bull. Soc. Math. France}, 144(3):395--427, 2016.

\bibitem{Goncharov:MTM}
A.~B. {Goncharov}.
\newblock {Multiple polylogarithms and mixed Tate motives}.
\newblock arXiv:math/0103059, 2001.

\bibitem{Hain:KZB}
R.~{Hain}.
\newblock {Notes on the Universal Elliptic KZB Equation}.
\newblock arXiv:1309.0580, 2013.

\bibitem{Hain:HodgeDeRham}
R.~Hain.
\newblock The {H}odge--de {R}ham theory of modular groups.
\newblock In {\em Recent advances in {H}odge theory}, volume 427 of {\em London
  Math. Soc. Lecture Note Ser.}, pages 422--514. Cambridge Univ. Press,
  Cambridge, 2016.

\bibitem{HM}
R.~{Hain} and M.~{Matsumoto}.
\newblock {Universal Mixed Elliptic Motives}.
\newblock arXiv:1512.03975, 2015.

\bibitem{IKZ}
K.~Ihara, M.~Kaneko, and D.~Zagier.
\newblock Derivation and double shuffle relations for multiple zeta values.
\newblock {\em Compos. Math.}, 142(2):307--338, 2006.

\bibitem{Ihara:Annals}
Y.~Ihara.
\newblock Profinite braid groups, {G}alois representations and complex
  multiplications.
\newblock {\em Ann. of Math. (2)}, 123(1):43--106, 1986.

\bibitem{Levin:Compositio}
A.~Levin.
\newblock Elliptic polylogarithms: an analytic theory.
\newblock {\em Compositio Math.}, 106(3):267--282, 1997.

\bibitem{LR}
A.~Levin and G.~Racinet.
\newblock Towards multiple elliptic polylogarithms.
\newblock arXiv:math/0703237, 2007.

\bibitem{LMS}
P.~Lochak, N.~Matthes, and L.~Schneps.
\newblock {Elliptic multiple zeta values and the elliptic double shuffle
  relations}.
\newblock arXiv:1703.09410, 2017.

\bibitem{Manin:Iterated}
Y.~I. Manin.
\newblock Iterated integrals of modular forms and noncommutative modular
  symbols.
\newblock In {\em Algebraic geometry and number theory}, volume 253 of {\em
  Progr. Math.}, pages 565--597. Birkh\"auser Boston, Boston, MA, 2006.

\bibitem{Matthes:Thesis}
N.~Matthes.
\newblock {\em {Elliptic multiple zeta values}}.
\newblock PhD thesis, Universit\"at Hamburg, 2016.

\bibitem{Matthes:Edzv}
N.~Matthes.
\newblock Elliptic double zeta values.
\newblock {\em J. Number Theory}, 171:227--251, 2017.

\bibitem{Nakamura:Galoisrep}
H.~Nakamura.
\newblock On exterior {G}alois representations associated with open elliptic
  curves.
\newblock {\em J. Math. Sci. Univ. Tokyo}, 2(1):197--231, 1995.

\bibitem{Nakamura:Eisenrevisited}
H.~Nakamura.
\newblock On profinite {E}isenstein periods in the monodromy of universal
  elliptic curves.
\newblock
  http://www.math.sci.osaka-u.ac.jp/$\sim$nakamura/zoo/fox/EisenRevisited.pdf,
  2016.

\bibitem{Reutenauer}
C.~Reutenauer.
\newblock {\em Free {L}ie algebras}, volume~7 of {\em London Mathematical
  Society Monographs. New Series}.
\newblock The Clarendon Press, Oxford University Press, New York, 1993.
\newblock Oxford Science Publications.

\bibitem{Serre:Lie}
J.-P. Serre.
\newblock {\em Lie algebras and {L}ie groups}, volume 1500 of {\em Lecture
  Notes in Mathematics}.
\newblock Springer-Verlag, Berlin, 2006.
\newblock 1964 lectures given at Harvard University, Corrected fifth printing
  of the second (1992) edition.

\bibitem{Tsunogai:Derivations}
H.~Tsunogai.
\newblock On some derivations of {L}ie algebras related to {G}alois
  representations.
\newblock {\em Publ. Res. Inst. Math. Sci.}, 31(1):113--134, 1995.

\bibitem{Zagier:Periods}
D.~Zagier.
\newblock Periods of modular forms and {J}acobi theta functions.
\newblock {\em Invent. Math.}, 104(3):449--465, 1991.

\bibitem{Zagier:Traces}
D.~Zagier.
\newblock Periods of modular forms, traces of {H}ecke operators, and multiple
  zeta values.
\newblock {\em S\=urikaisekikenky\=usho K\=oky\=uroku}, (843):162--170, 1993.
\newblock Research into automorphic forms and $L$ functions (Japanese) (Kyoto,
  1992).

\end{thebibliography}

\end{document}